\documentclass{svjour3}                     
\smartqed  
\usepackage{amsmath}
\usepackage{amssymb}
\usepackage{latexsym}
\usepackage{graphicx}
\usepackage{mathptmx}      
\usepackage{bm}           
\newcommand{\R}{\,{\mathbb R}}
\newcommand{\diag}{\mbox{diag}}

\newcommand{\sign}{\mbox{sign}}

\newcommand{\bolda}{{\mathbf a}}
\newcommand{\bolde}{{\mathbf e}}

\newcommand{\boldq}{{\mathbf q}}
\newcommand{\boldu}{{\mathbf u}}
\newcommand{\boldx}{{\mathbf x}}

\newcommand{\boldG}{{\mathbf G}}
\newcommand{\boldH}{{\mathbf H}}

\newcommand{\boldphi}{\boldsymbol{\varphi}}
\newcommand{\boldpsi}{\boldsymbol{\psi}}
\newcommand{\boldzero}{{\mathbf 0}}
\newcommand{\binomial}[2]{\begin{pmatrix}#1\\#2\end{pmatrix}}

\journalname{Multidimensional Systems and Signal Processing}
\begin{document}

\title{Matrix Spectral Factorization - SA4 Multiwavelet%
}


\author{Vasil Kolev \and Todor Cooklev \and Fritz~Keinert}


\institute{Vasil Kolev \at
  Institute of Information and Communication Technologies \\
  Bulgarian Academy of Sciences \\
  Bl. 2 Acad. G. Bonchev St. \\
  1113 Sofia, Bulgaria \\
  \email{kolev\_acad@abv.bg} \\
  \and
  Todor Cooklev, {\em Senior Member IEEE} \at
  Wireless Technology Center \\
  Indiana University -- Purdue University \\
  Fort Wayne, IN 46805, USA \\
  \email{cooklevt@ipfw.edu} \\
  \and
  Fritz Keinert \at
  Dept.~of Mathematics \\
  Iowa State University \\
  Ames, IA 50011, USA \\
  \email{keinert@iastate.edu}
}

\date{Received: date / Accepted: date}

\maketitle

\begin{abstract}
  In this paper, we investigate Bauer's method for the matrix spectral
  factorization of an $r$-channel matrix product filter which is a
  half-band autocorrelation matrix. We regularize the resulting matrix
  spectral factors by an averaging approach and by multiplication by a
  unitary matrix. This leads to the approximate and exact orthogonal
  SA4 multiscaling functions. We also find the corresponding
  orthogonal multiwavelet functions, based on the QR decomposition.

\end{abstract}

\section{Introduction}\label{sec:intro}

If $H(z) = C_0 + C_1 z^{-1} + \cdots C_n z^{-n}$ is a matrix
polynomial which represents a causal multifilter bank, the product
filter bank is given by $P(z) = H(z) H^*(z)$, where $H^*(z) =
H^T(z^{-1})$.  The coefficients $C_k$ are $r \times r$ matrices.
{\em Matrix Spectral Factorization} (MSF) describes the problem of
finding $H(z)$, given $P(z)$. 

MSF is still quite unknown in the signal processing community. Because
of its numerous possible applications, in particular in multiple-input
and multiple-output (MIMO) communications, image processing,
multidimensional control theory, and others, it deserves to receive
greater attention. We are interested in MSF as a tool for designing
orthogonal multiwavelets.

It is known that MSF is possible as long as $P(z)$ is positive
semidefinite for $z$ in the unit circle in the complex plane
\cite{EJL-09,HHS-04}.  We call $P(z)$ {\em degenerate} if it is not
strictly positive definite, that is, $\det(P(z)) = 0$ at one or more
points.

A number of numerical approaches to MSF have been proposed, but they
usually cannot handle the degenerate case.  Unfortunately, those are
exactly the cases of greatest interest for applications in the
construction of multiwavelets.

We use Bauer's method, as described by Youla and Kazanjian
\cite{YK-78}, which is based on the Cholesky factorization of a block
Toeplitz matrix. This algorithm can handle the degenerate case, but
convergence is quite slow.

We study the SA4 multiwavelet as a test case, and as a benchmark for
future applications. Starting from a known $H(z)$, we compute $P(z)$
and factor it again.

It is easy to see that matrix spectral factorization is not unique. If
$H_1(z)$ is any factor, then $H_2(z) = H_1(z) U$ is also a factor for
any orthogonal matrix $U$, or more generally $H_2(z) = H_1(z) U(z)$
for paraunitary $U(z)$.  For the factor produced by Bauer's method,
the matrix $C_0$ is always lower triangular, so we are forced to use
further modifications to recover the original SA4.

In this paper, we only concentrate on modifications which speed up
convergence and produce factors close to original SA4. This leads to
two factorizations: the approximate SA4 multifilter, which is close to
SA4, and the exact SA4 multifilter, which is equal to SA4 except for
insignificant roundoff. In future papers, we plan to impose desirable
conditions directly on the spectral factor, without working towards a
known result.

The main contributions of this paper are

\begin{itemize}
  \item Considering a novel method for computing the multichannel
    spectral factorization of degenerate matrix polynomial matrices;
  \item Speeding up the slow convergence of the algorithm;
  \item Minimizing the numerical errors which appear in the algorithm;
  \item Demonstrating by example that this approach can be used to
    construct an orthogonal symmetric multiwavelet filter;
  \item Comparing the frequency responses and experimental numerical
    performance of the approximate and exact SA4 multiwavelets.
\end{itemize}

The rest of the paper is organized as follows.  In
chapter~\ref{sec:background} we give some background material, state
the problem, and describe related research. We also introduce matrix
product filter theory, and present the product filter of the SA4
multiscaling function as a benchmark test case. In
chapter~\ref{sec:spectfact} we consider Bauer's method, and discuss
its numerical behavior. In chapter~\ref{sec:construction} we describe
the approach for obtaining the approximate and exact SA4 multiscaling
and multiwavelet functions. The performance analysis of these
multiwavelets is shown in chapter~\ref{sec:perf}.
Chapter~\ref{sec:concl} gives conclusions.

\section{Background and Problem Statement}\label{sec:background}

We will use the following notation conventions, illustrated with the letter ‘a': 

$a$ -- lowercase letters refer to scalars or scalar-valued functions; 

$\bolda$ -- lowercase bold letters are vectors or vector-valued functions; 

$A$ -- uppercase letters are matrices or matrix-valued functions. 

\noindent The symbols $I$ and $0$ denote the identity and zero matrices of
appropriate size, respectively.

We will be using polynomials in a complex variable $z$ on the unit
circle, but all coefficients will be real-valued. Thus, the complex
conjugate $A^*(z)$ of a matrix polynomial $A(z) = \sum_k A_k z^{k}$ is
given by
\begin{equation*}
  A^*(z) = A^T(z^{-1}) = \sum_k A_k^T z^{-k}.
\end{equation*}

\subsection{Matrix Product Filters}

\begin{figure}
  \includegraphics[width=4in]{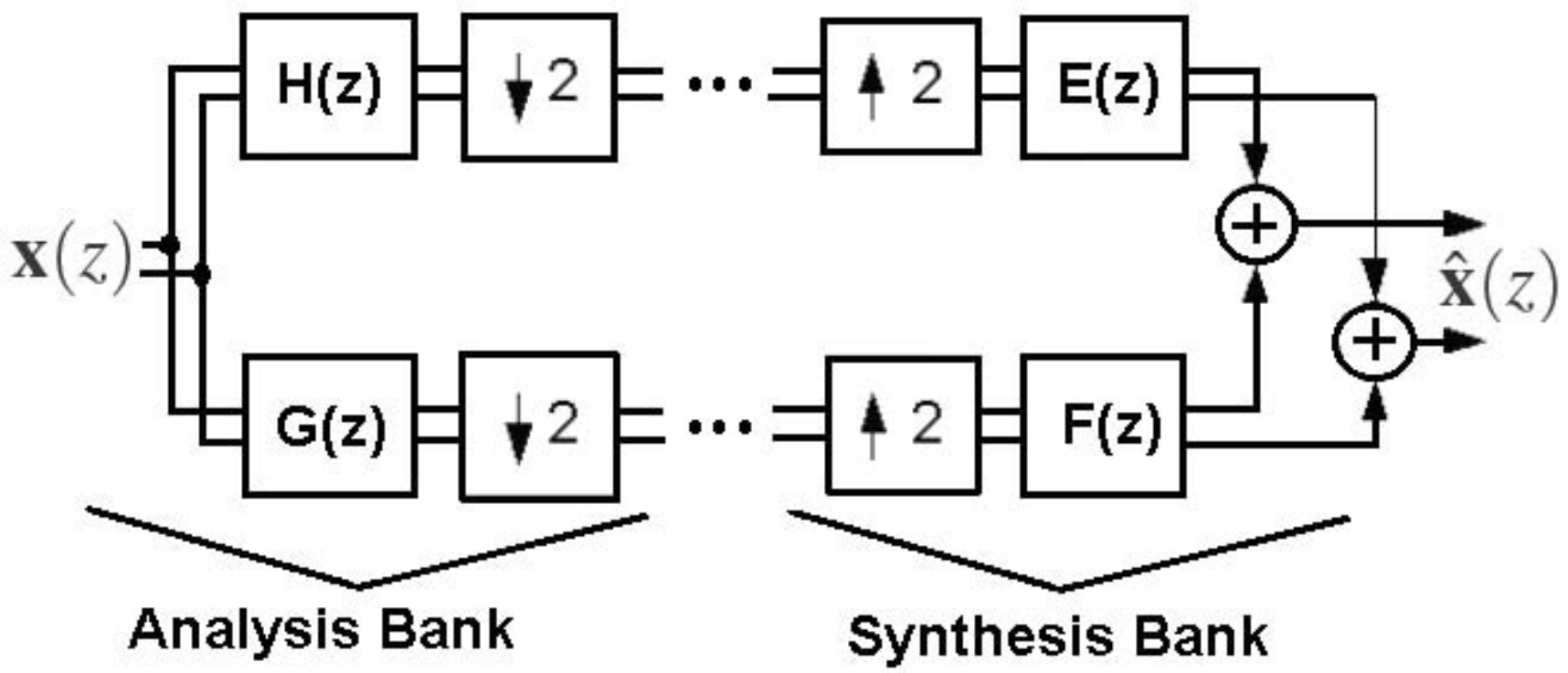} \\
  \caption{Two-channel multifilter bank.}
  \label{fig:01}
\end{figure}

A two-channel multiple-input multiple-output (MIMO) filter bank of
multiplicity $r$ is shown in fig.~\ref{fig:01}. Here $\boldx(z)$ is
the input signal vector of length $r$, and the analysis multifilters
are
\begin{align*}
  H(z) &= C_0 + C_1 z^{-1} + \cdots C_n z^{-n}, \\
  G(z) &= D_0 + D_1 z^{-1} + \cdots D_n z^{-n},
\end{align*}
where $C_k$, $D_k$ are matrices of size $r \times r$. Likewise, $E(z)$
and $F(z)$ are the synthesis multifilters, and $\hat\boldx(z)$ is the
output vector signal. For simplicity, we assume that all coefficients
are real.

The input-output relation of this filter bank is
\begin{equation*}
  \hat\boldx(z) = \frac{1}{2} \left[ E(z) H(z) + F(z) G(z) \right]
  \boldx(z) + \frac{1}{2} \left[ E(z) H(-z) + F(z) G(-z) \right] \boldx(-z).
\end{equation*}
It is convenient to write this equation in matrix form as
\begin{equation*}
  \begin{bmatrix}
    \hat\boldx(z) \\ \hat\boldx(-z)
  \end{bmatrix} = \frac{1}{2}
  \begin{bmatrix}
    E(z) & F(z) \\
    E(-z) & F(-z)
  \end{bmatrix}
  \begin{bmatrix}
    H(z) & H(-z) \\
    G(z) & G(-z)
  \end{bmatrix}
  \begin{bmatrix}
    \boldx(z) \\ \boldx(-z)
  \end{bmatrix}.
\end{equation*}

In general, the design of the multifilter bank requires four
multifilters: two on the analysis and two on the synthesis side. In
perfect-reconstruction orthogonal MIMO filter banks, the analysis
modulation matrix
\begin{equation*}
  M(z) =
  \begin{bmatrix}
    H(z) & H(-z) \\
    G(z) & G(-z)
  \end{bmatrix}.
\end{equation*}
is {\em paraunitary}, that is
\begin{equation}\label{eq:03}
  M^*(z) M(z) = M(z) M^*(z) = I.
\end{equation}
In this case we can define the synthesis filters in terms of the
analysis filters:
\begin{equation*}
  E(z) = H^*(z), \quad F(z) = G^*(z).
\end{equation*}
Eq.~\eqref{eq:03} can also be expressed as
\begin{equation*}
   \begin{split}
     H(z) H(z)^* + H(-z) H(-z)^* = I, \\
     G(z) G(z)^* + G(-z) G(-z)^* = I, \\
     H(z) G(z)^* + H(-z) G(-z)^* = 0,
   \end{split}
\end{equation*}
or in terms of the coefficients
\begin{equation*}
  \begin{split}
    \sum_k C_k C_{k+2\ell}^T &= \sum_k D_k D_{k+2\ell}^T = \delta_{0\ell} \,I, \\
    \sum_k C_k D_{k+2\ell}^T &= 0,
  \end{split}
\end{equation*}
for any integer $\ell$.

Given $H(z)$, the matrix polynomial
\begin{equation*}
  P(z) = H(z) H^*(z) = \left( \sum_{k=0}^n C_k z^{-k} \right)
  \left( \sum_{k=0}^n C_k^T z^k \right) = \sum_{k=-n}^n P_k z^k,
\end{equation*}
is called the {\em matrix lowpass product filter} (or scalar product
filter in the case $r=1$).  The coefficients $P_k$ satisfy $P_{-k} =
P_k^T$. $P(z)$ is a half-band polynomial filter \cite{CNKS-96}, that
is, it satisfies the equation
\begin{equation}\label{eq:08}
  P(z) + P(-z) = 2I.
\end{equation}
This implies that $P_0 = I$, and $P_{2\ell} = 0$ for all nonzero
integers $\ell$.

\subsection{Multiwavelet Theory}\label{subsec:mwtheory}

\subsubsection{Multiscaling and Multiwavelet Functions}

Consider the iteration of a MIMO filter bank along the channel
containing $H(z)$. After $k$ iterations, the equivalent filters will
be
\begin{equation*}
   \begin{split}
     H^{(k)}(z) &= \prod_{j=0}^k H(z^{2^j}) = H(z^{2^k}) H(z^{2^{k-1}}) \cdots H(z), \\
     G^{(k)}(z) &= G(z^{2^k}) \prod_{j=0}^{k-1} H(z^{2^j}) = G(z^{2^k}) H(z^{2^{k-1}}) \cdots H(z),
   \end{split}
\end{equation*}
or in the time domain
\begin{align*}
  C_j^{(k)} &= \sum_j C_j^{(k-1)} C_{j-2^{k-1}}^{(k-1)}, \\
  D_j^{(k)} &= \sum_j D_j^{(k-1)} D_{j-2^{k-1}}^{(k-1)},
\end{align*}
where $C_j^{(0)} = C_j$, $D_j^{(0)} = D_j$. 

The filterbank coefficients are associated with function vectors
$\boldphi = [\phi_0, \phi_1, \ldots, \phi_{r-1}]^T$, called the {\em
  multiscaling function}, and $\boldpsi = [\psi_0, \psi_1, \ldots,
  \psi_{r-1}]^T$, called the {\em multiwavelet function}. These
functions satisfy the recursion equations
\begin{equation*}
  \begin{split}
    \boldphi(t) &= \sqrt{2} \sum_{k=0}^n C_k \boldphi(2t-k), \\
    \boldpsi(t) &= \sqrt{2} \sum_{k=0}^n D_k \boldphi(2t-k).
  \end{split}
\end{equation*}
The support of $\boldphi$, $\boldpsi$ is contained in the interval
$[0,n]$, but it could be strictly smaller than the interval.
See \cite{MRvF-96} for more details.

\subsubsection{Multiwavelet Properties}

\hspace{\parindent}{\bf Pre- and Postfiltering}

Since multifilters are MIMO, they operate on several streams of input
data rather than one. Therefore, input data needs to be vectorized.
There are many prefilters available for various applications, but
three kinds are preferred:
\begin{itemize}
  \item Oversampling - Repeated Row Preprocessing 
  \item Critical Sampling  - Matrix Preprocessing
  \item Embedded Orthogonal Symmetric Prefilter Bank
\end{itemize}

In \cite{SW-00}, the first and second approach is presented for the
GHM and CL multiwavelets. In practical applications, a popular choice
is based on the Haar transform \cite{CMP-98,TST-99}.  Haar pre- and
postfilters have the advantage of simultaneously possessing symmetry
and the orthogonality, and no multiplication is needed if one ignores
the scaling factor. The problem of constructing symmetric prefilters
is considered in detail in \cite{HSLS-03}, with proposed solutions for
the DGHM and CL multifilters.

The third technique is applied in \cite{GM-05,HLSH-07}.  They find a
three-tap prefilter for the DGHM multifilter by searching among the
parameters which minimize the first-order absolute moment of the filter
coefficients.

Similarly, at the filter bank output, a postfilter is needed to
reconstruct the signal. Pre- and postprocessing is not needed for
scalar wavelets. \\

{\bf Balancing Order}

The orthogonal multiscaling function is said to be {\em balanced of
  order $q$} if the signals $\boldu_k = [\cdots, (-2)^k, (-1)^k, 0^k,
  1^k, 2^k, \cdots]^T$ are preserved by the operator $L^T$ for $k =
0,1,\ldots,q-1$. That is,
\begin{equation*}
  L^T \boldu_k = 2^{-k} \boldu_k,
\end{equation*}
where
\begin{equation*}
  L =
  \begin{bmatrix}
    \cdots \\
    & C_0 & C_1 & \cdots & C_n \\
    & & C_0 & C_1 & \cdots & C_n \\
    & & & & & & \cdots
  \end{bmatrix}.
\end{equation*}
Multiwavelets that do not satisfy this property are said to be {\em
  unbalanced}. See \cite{LV-01,LV-98,PS-98,LP-11} for more
details. Balanced multiwavelets do not require preprocessing. \\

{\bf Approximation Order}

In the scalar case, a certain approximation order refers to the
ability of the low-pass filter to reproduce discrete-time polynomials,
while the wavelet annihilates discrete-time polynomials up to the same
degree.  Since many real-world signals can be modeled as polynomials,
this property typically leads to higher coding gain (CG).

In the case of non-balanced multiwavelets, appropriate preprocessing
is required to take advantage of approximation orders.

 The approximation order of multiscaling and multiwavelet functions
 can be determined using the following result established in
 \cite{WLXL-10}. A multiscaling function provides approximation order
 $p$ iff there exist vectors $\boldu_k \in \R^r$, $0 \le k < p$,
 $\boldu_0 \ne \boldzero$, which satisfy
\begin{equation*}
   \begin{split}
     \sum_{\ell=0}^k \frac{1}{(2i)^\ell} \binomial{k}{\ell}
     \boldu_{k-\ell}^T D^\ell H(0) &= 2^{-k} \boldu_k^T, \\
     \sum_{\ell=0}^k \frac{1}{(2i)^\ell} \binomial{k}{\ell}
     \boldu_{k-\ell}^T D^\ell H(\pi) &= 0,
   \end{split}
\end{equation*}
where $D^\ell$ is the derivative of order $\ell$, and the
{\em masks} of multiscaling and multiwavelet functions are
\begin{equation}\label{eq:mask}
  H(\omega) = \frac{1}{\sqrt{2}} \sum_{k=0}^n C_k e^{i k \omega}, \qquad
  G(\omega) = \frac{1}{\sqrt{2}} \sum_{k=0}^n D_k e^{i k \omega}.
\end{equation}
To fully characterize the multifilter bank with respect to
approximation order we can add the condition
\begin{equation*}
  \sum_{\ell=0}^k \frac{1}{(2i)^\ell} \binomial{k}{\ell}
  \boldu_{k-\ell}^T D^\ell G(0) = 0.
\end{equation*}\\

{\bf Good Multifilter Properties (GMPs)}

A multiwavelet system can be represented by an equivalent scalar
filter bank system \cite{TSLT-00}. The $r$ equivalent scalar filters
are, in fact, the $r$ polyphases of the corresponding multifilter.
This relationship motivates a new measure called {\em good multifilter
  properties} (GMPs) \cite{TSLT-98}, \cite[(Def.~1)]{TSLT-00}, which
characterizes the magnitude response of the equivalent scalar filter
bank associated with a multifilter. GMPs provide a set of design
criteria imposed on the scalar filters, which can be translated
directly to eigenvector properties of the designed multiwavelet
filters.

An orthogonal multiwavelet system has GMPs of order at least $(1,1,1)$
if the following conditions are satisfied \cite{TSLT-00}
\begin{equation*}
   \begin{split}
     \boldH(0) \bolde(0) = \bolde(0) \\
     \boldH(r \pi) \bolde(\pi)=0 \\
     \boldG(0) \bolde(0)=0 
   \end{split}
\end{equation*}
where $H(\omega)$ and $G(\omega)$ are masks of the lowpass and
highpass filters (see eq.~\eqref{eq:mask}), and
\begin{equation*}
  \bolde(\omega) = [1, e^{-i\omega} , e^{-2i\omega}, \ldots, e^{-(r-1)i\omega} ]^T.
\end{equation*}
A class of symmetric-antisymmetric biorthogonal multiwavelet filters
which possess GMPs was introduced in \cite{TST-99}.

A GMP order of at least $(1,1,1)$ is critical for ensuring no
frequency leakage across bands, hence improving compression
performance.  Multifilters possessing GMPs have better performance
than those which do not possess GMPs \cite{LY-10}, due to better
frequency responses for energy compaction, greater regularity and
greater approximation order of the corresponding wavelet/scaling
functions. GMPs help prevent both DC and high-frequency leakage across
bands; this contributes to reduced smearing, blocking, ringing
artifacts, and also helps to prevent checkerboard artifacts in
reconstructed images for image coding.

The SA4 multiwavelet, introduced below in section
\ref{subsec:testing}, is a member of a one-parameter family of
orthogonal multiwavelets with four coefficient matrices
\cite{TSLT-00}. Different members of the family have different
approximation order, but all of them have a GMP order of at least $(1,
1, 1)$.  This manifests itself in the smooth decay to
zero of the magnitude responses near $\omega = \pi$, compared to the
following 4-tap orthogonal multiwavelet filters:

\begin{itemize}
  \item The GHM multiwavelet \cite{DGHM-96} has symmetric orthogonal
    scaling functions and an approximation order of 2 for its filter
    length 4.
  \item Chui and Lian's CL multiwavelet \cite{CL-96} has the highest
    possible approximation order of 3 for its filter length 3.
  \item Jiang's multiwavelet JOPT4 \cite{J-98a} has optimal
    time-frequency localization for its filter length. 
\end{itemize}

Although the GHM and CL systems are the most commonly used orthogonal
multiwavelet systems, and have higher approximation order than the SA4
system, they do not satisfy GMPs.\\

{\bf Smoothness}

The smoothness of $\boldphi$, $\boldpsi$ can be characterized by
Sobolev regularity $S$; this is discussed in
\cite{MS-97,CDP-97,J-98b}.  The regularity estimate is related to both
the eigenvalues and the corresponding eigenvectors of the transition
operator and to spectral radius of the transition operator.\\

{\bf Symmetry and Antisymmetry}

A function $g(t)$ is {\em symmetric} about a point $c$ if $g(c-t) =
g(c+t)$ for all $t$. It is {\em antisymmetric} if $g(c-t) = - g(c+t)$.

For scaling functions and wavelets of compact support, the only
possible symmetry is about the center of their support. For
multiscaling and multiwavelet functions, this could be a different
point for different components.

In the scalar case, the scaling function cannot be antisymmetric,
since it must have a nonzero integral. In the multiwavelet case, some
components even of $\boldphi$ can have antisymmetry.

\subsection{Problem Statement}

The design of perfect reconstruction multifilter banks and
multiwavelets remains a significant problem. In the scalar case,
spectral factorization of a half-band filter that is positive definite
on the unit circle was, in fact, the first design technique, suggested
by Smith and Barnwell \cite{SB-86}. This provides a motivation to
extend the technique to the multiwavelet case.  Spectral factorization
of a Hermitian product filter in the multiwavelet case is more
challenging than in the scalar case, and has not been investigated
previously.

The present paper considers the application of MSF to the problem of
finding an orthogonal lowpass multifilter $H(z)$, given a product
filter $P(z)$.  Since every orthogonal multiscaling filter $H(z)$ is a
spectral factor of some matrix product filter, this can be used as a
tool for designing new filters.

With an eye towards future applications, we are interested in
constructing multifilters with desirable properties, especially
GMPs. This implies that the determinant of $P(z)$ will have at least
one zero of higher multiplicity on the complex unit circle.  That is,
$P(z)$ will be degenerate.

As a test case, we want to use a product filter $P(z)$ which satisfies
the half-band condition \eqref{eq:08}, and is derived from an $H(z)$
which has more than 2 taps (i.e. $n \ge 3$) and good regularity
properties and GMPs. Our benchmark test case will be the 2-channel
orthogonal SA4 multiwavelet.

While there are many matrix spectral factorization algorithms
\cite{CNKS-96,JK-85,L-93}, most cannot handle the degenerate case. We
use Bauer's method, based on the Toeplitz method of spectral
factorization of the Youla and Kazanjian algorithm \cite{B-56,YK-78},
which can handle this case.  However, convergence becomes very slow.

It is easy to see that matrix spectral factorization is not unique. If
$H_1(z)$ is a factor of a given $P(z)$, then $H_2(z) = H_1(z) U$ is
also a factor for any orthogonal matrix $U$, or more generally $H_2(z)
= H_1(z) U(z)$ for paraunitary $U(z)$. Other operations may be
permissible in some cases.

In our experiments with Bauer's algorithm, the initial factors do not
correspond to smooth functions. We apply regularization techniques,
both to speed up convergence and to work towards recovering the
original filter in this test case. In this manner, we derive two
filters which we call the {\em approximate} and {\em exact SA4
  multiwavelets}.

Let us briefly summarize some previous work.

Matrix Spectral Factorization (MSF) plays a crucial role in the
solution of various applied problems for MIMO systems in
communications and control engineering \cite{WMW-15}. It has been
applied to designing minimum phase FIR filters and the associated
all-pass filter \cite{HCW-10}, quadrature-mirror filter banks
\cite{C-07}, MIMO systems for optimum transmission and reception
filter matrices for precoding and equalization \cite{F-05}, precoders
\cite{DXD-13}, and many other applications.

The analysis of linear systems corresponding to a given spectral
density function was first established by Wiener, who used linear
prediction theory of multi-dimensional stochastic processes
\cite{WM-57}. The method of Wilson for MSF was developed for
applications in the prediction theory of multivariate discrete time
series, with known (or estimated) spectral density \cite{W-72}. Many
numerical approaches to MSF have been proposed, for example
\cite{CNKS-96,JK-85,L-93}. A survey of such algorithms is given in
\cite{SK-01}; however, this does not include algorithms for the
degenerate case.

It is well known that all MSF methods have difficulties in this
situation, and some of them cannot handle zeros on the unit circle at
all.  For example, Ku{\v c}era's algorithm, an otherwise popular matrix
spectral factorization algorithm, has this limitation, and is not
suitable for our purpose.

Bauer's method was successful applied in \cite{LR-86} to the Radon
projection. In that paper, factorization of the autocorrelation matrix
of the Radon projection of a minimum phase pseudo-polynomial restored
the coefficients of the original pseudo-polynomial with an accuracy
around $10^{-17}$ after 10,000 recursions on a 64-bit floating point
processor.

\section{Matrix Spectral Factorization}\label{sec:spectfact}

The following theorem answers the existence question for MSF.

\begin{theorem} (Matrix form of the Riesz-Fej{\'e}r lemma \cite{HHS-04,EJL-09})
  A matrix polynomial
  \begin{equation*}
    P(z) = \sum_{k=-n}^n P_k z^k
  \end{equation*}
  can be factored as  $P(z) = H(z) H^*(z)$, where
  \begin{equation*}
    H(z) = \sum_{k=0}^n C_k z^{-k}
  \end{equation*}
  is an $r$-channel causal polynomial matrix, if and only if $P(z)$ is
  symmetric positive semidefinite for $z$ on the complex unit circle.
\end{theorem}
	
\subsection{Bauer's Method}

If $H(z)$ and $P(z)$ are known, we can construct doubly infinite block
Toeplitz matrices $L$ and $T$ from their coefficients, by setting
$L_{ij} = C_{i-j}$, $T_{ij} = T_{j-i}$. $T$ is symmetric and block
banded with bandwidth $n$; $L$ is block lower triangular, also with
bandwidth $n$.

The relation $P(z) = H(z) H^*(z)$ corresponds to $T = L \ L^T$, that
is, $L$ is a Cholesky factor of $T$.

In Bauer's method, we pick a large enough integer $f$, and truncate
$T$ to (block) size $f \times f$:
\begin{equation*}
  T^{(f)} =
  \begin{pmatrix}
    P_0    & P_1   & \cdots  & P_n & \\
    P_{-1} & P_0   & P_1   & \cdots   & P_n  & \\
    \vdots & & \ddots & & & \ddots &  \\
    P_{-n} &  & & \ddots & &  & P_n \\
          & \ddots & & & \ddots & & \vdots \\
     &  & \ddots & &  & \ddots & P_1 \\
     &  &  & P_{-n} & \cdots & P_{-1} & P_0
  \end{pmatrix}
\end{equation*}

If $T^{(f)}$ is positive definite, we can compute the Cholesky
factorization $T = L^{(f)} \ L^{(f)T}$, where
\begin{equation*}
  L^{(f)} =
  \begin{pmatrix}
    C_0^{(1)}  & \\
    C_1^{(2)} & C_0^{(2)} & \\
    \vdots & & \ddots &  & &  & \\
    C_n^{(n+1)} &  & & C_0^{(n+1)} & \\
     & \ddots & &  & \ddots &  \\
     &  & C_n^{(f)} & \cdots & C_1^{(f)}  &C_0^{(f)}
  \end{pmatrix}
\end{equation*}

It is shown in \cite{YK-78} that this factorization is possible, and
that $C_k^{(f)} \to C_k$ as $f \to \infty$.

Bauer's method works even in the case of highly degenerate $P(z)$.
Unfortunately, convergence in this case is very slow.

\subsection{Numerical Behavior}\label{subsec:numerical}

\begin{figure}
  \includegraphics[width=2.5in]{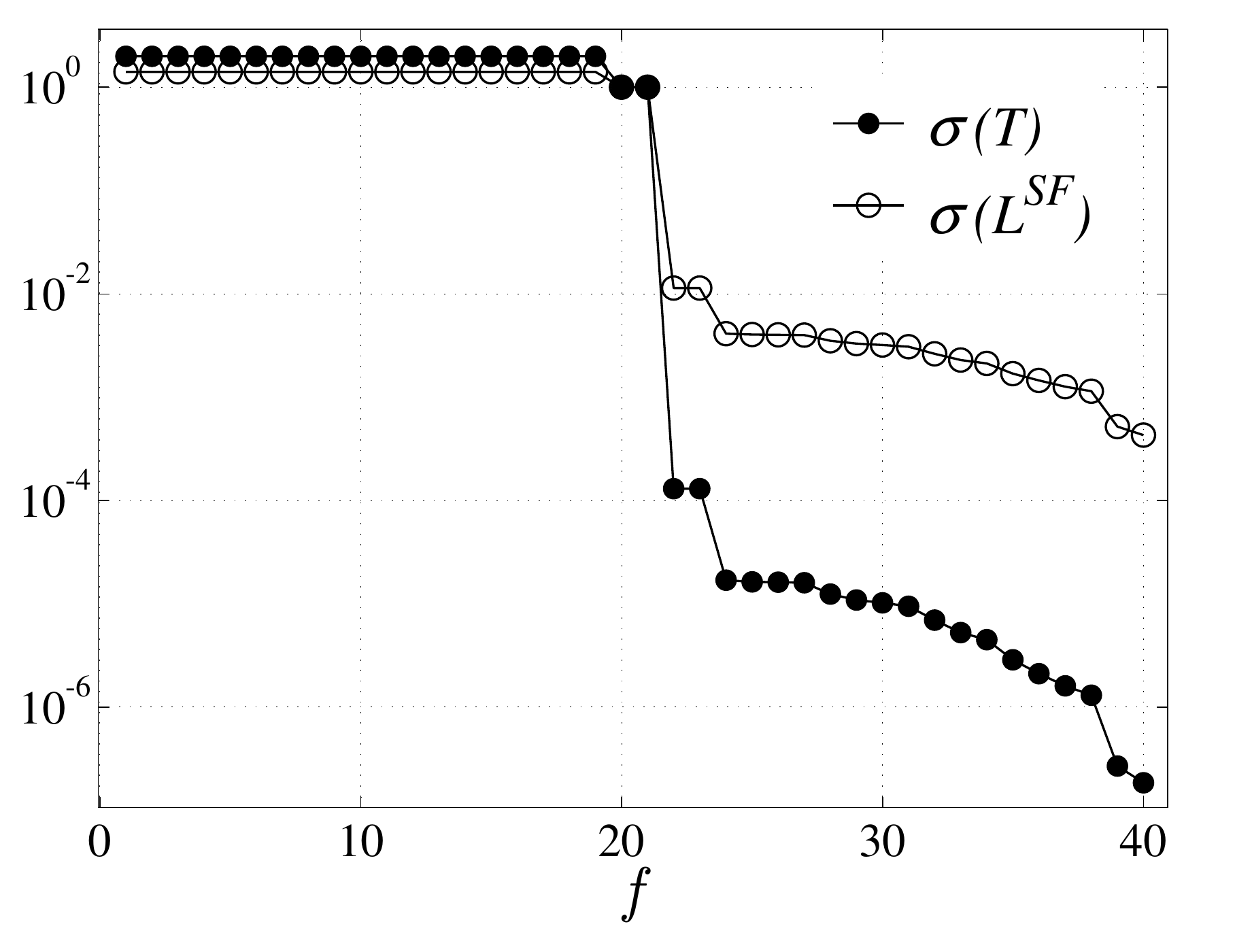}
  \caption{A sharp drop of the singular values of the Toeplitz matrix
    $T$ and of its spectral factor $L^{SF}$, for $f=20$; (filled
    circles) singular values of $T$; (hollow circles) singular values
    of $L^{SF}$.}
  \label{fig:02}
\end{figure}

For chosen size $f$, the matrix $T^{(f)}$ is of size $2f \times 2f$.  Its
singular values are all in the range $[0,2]$. The first $(f-1)$
singular values are close to 2, then $\sigma_f = \sigma_{f+1} = 1$,
and the remaining $(f-1)$ singular values are close to 0.  See
fig.~\ref{fig:02}.

The numerical stability of the factorization algorithm is good, until
the higher singular values get too close to 0.

A bigger problem is the slow convergence. For example, in the scalar
case $p(z) = z^{-1} + 2 + z$, where the factor is $h(z) = 1 + z^{-1}$,
the value on the diagonal in row $f$ is
\begin{equation*}
  \sqrt{1 + \frac{1}{f}} \approx 1 + \frac{1}{2f},
\end{equation*}
which converges only very slowly to the limit 1.  This simple problem
has a double zero of the determinant on the unit circle. For higher
order zeros, convergence speed becomes even worse.

\subsection{Testing Bauer's Method}\label{subsec:testing}

\begin{figure}
  \includegraphics[width=3.5in]{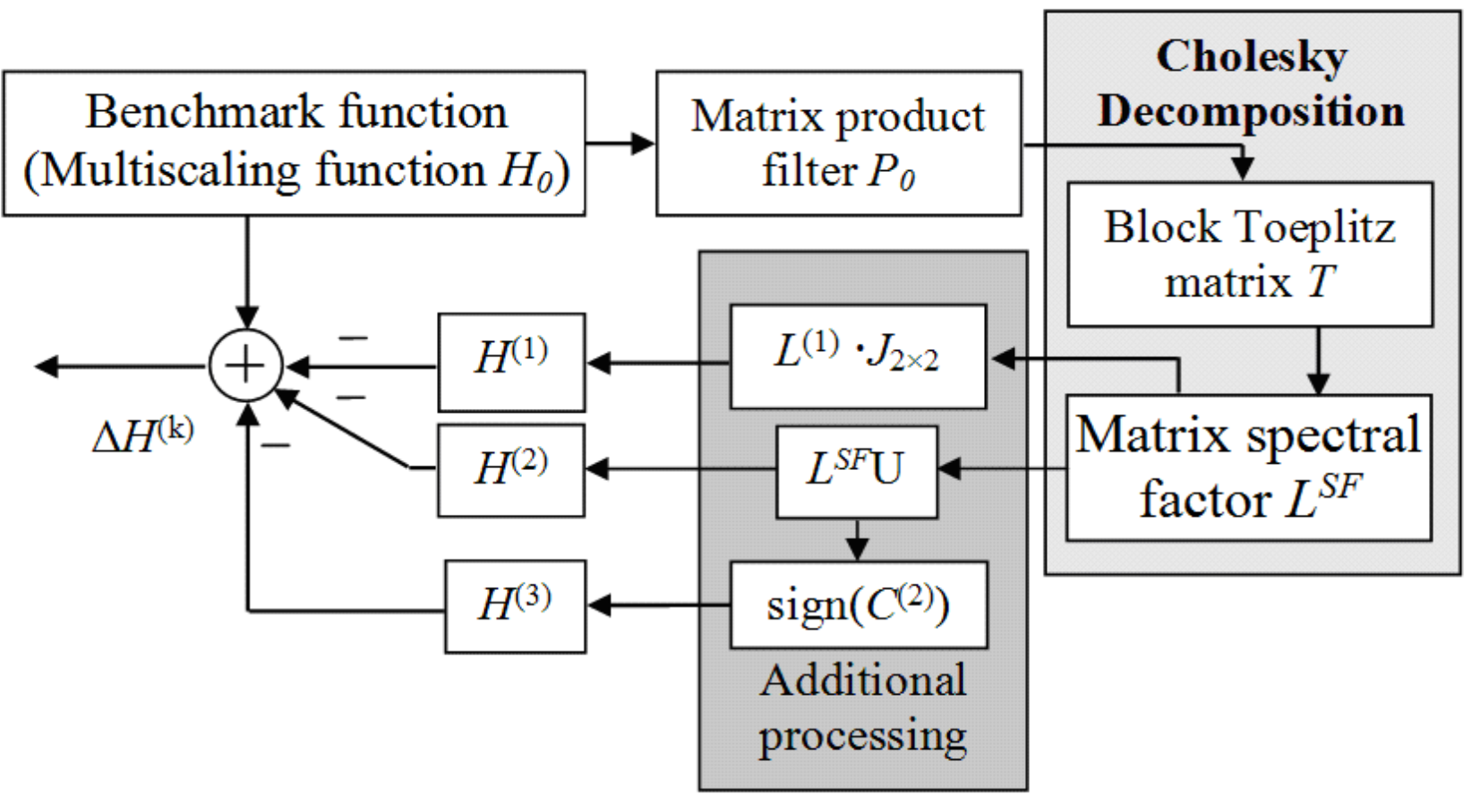} 
  \caption{Steps in testing matrix spectral factorization using Bauer's method.}
  \label{fig:03}
\end{figure}

A flowchart for testing Bauer's method for MSF is shown in
fig.~\ref{fig:03}.  The steps are

\begin{itemize}
  \item[Step 1:] Choice of the benchmark multiscaling function $H(z)$;
  \item[Step 2:] Construction of the matrix lowpass product filter $P(z)$; 
  \item[Step 3:] Find a spectral factor $L^{SF} = L^{(f)}$, and assess
    its quality;
  \item[Step 4:] Do different kinds of postprocessing to obtain spectral
    factors $H^{(i)}$, $i=1,2,3$;
  \item[Step 5:] Compare the obtained factors with the benchmark $H(z)$.
\end{itemize}

{\bf Step 1:} We consider the well-known orthogonal SA4 multiwavelet
\cite{CL-96,HM-12,TSLT-00} as a benchmark testing case.

The recursion coefficients $C_k$ and $D_k$ of SA4 depend on a
parameter $t$. They are
\begin{equation}\label{eq:24}
    \begin{gathered}
    C_0 = \alpha
    \begin{pmatrix}
      1 & t \\
      1 & -t
    \end{pmatrix}; \quad
    C_1 = \alpha
    \begin{pmatrix}
      t^2 & t \\
      -t^2 & t
    \end{pmatrix}; \quad
    C_2 = \alpha
    \begin{pmatrix}
      t^2 & -t \\
      t^2 & t
    \end{pmatrix}; \quad
    C_3 = \alpha
    \begin{pmatrix}
      1 & -t \\
      -1 & -t
    \end{pmatrix}; \\
    D_0 = \alpha
    \begin{pmatrix}
      t & -1 \\
      t & 1
    \end{pmatrix}; \quad
    D_1 = \alpha
    \begin{pmatrix}
      -t & t^2 \\
      t & t^2
    \end{pmatrix}; \quad
    D_2 = \alpha
    \begin{pmatrix}
      -t & -t^2 \\
      -t & t^2
    \end{pmatrix}; \quad
    D_0 = \alpha
    \begin{pmatrix}
      t & 1 \\
      -t & 1
    \end{pmatrix};
  \end{gathered}
\end{equation}
where $\alpha = 1/(\sqrt{2} (1+t^2))$. In this paper, we use the
filter with $t = 4 + \sqrt{15}$, which leads to $\alpha =
(4-\sqrt{15})/(8 \sqrt{2})$. Numerical values of $C_k$, $D_k$ are
listed in table \ref{tab:01}.

These coefficients generate the 2-band compactly supported orthogonal
multiscaling function $\boldphi = [\phi_0, \phi_1]^T$ and multiwavelet
function $\boldpsi = [\psi_0, \psi_1]^T$, all with support $[0,3]$.
Graphs of these functions can be found in fig.~\ref{fig:08} in a later
chapter.

\begin{table}
  \caption{Coefficients of the original SA4 multiscaling and multiwavelet functions.}
  \label{tab:01}
  \begin{tabular}{|c|rr|rr|}
    \hline
    $n$ & \multicolumn{2}{c|}{$C_n$} & \multicolumn{2}{c|}{$D_n$} \\
    \hline
    $0$ 
    & $0.011226792152545$ &  $0.088388347648318$ & $-0.088388347648318$ &  $0.011226792152545$ \\
    & $0.011226792152545$ & $-0.088388347648318$ & $-0.088388347648318$ & $-0.011226792152545$ \\
    \hline
    $1$
    &  $0.695879989034003$ & $0.088388347648318$ &  $0.088388347648318$ & $-0.695879989034003$ \\
    & $-0.695879989034003$ & $0.088388347648318$ & $-0.088388347648318$ & $-0.695879989034003$ \\
    \hline
    $2$
    & $0.695879989034003$ & $-0.088388347648318$ & $0.088388347648318$ &  $0.695879989034003$ \\
    & $0.695879989034003$ &  $0.088388347648318$ & $0.088388347648318$ & $-0.695879989034003$ \\
    \hline
    $3$
    &  $0.011226792152545$ & $-0.088388347648318$ & $-0.088388347648318$ & $-0.011226792152545$ \\
    & $-0.011226792152545$ & $-0.088388347648318$ &  $0.088388347648318$ & $-0.011226792152545$ \\
    \hline
  \end{tabular}
\end{table}

{\bf Step 2:} The symbol of the product lowpass multifilter $P$ has
the form
\begin{equation*}
  P(z) = P_3^T z^{-3} + P_1^T z^{-1} + P_0 + P_1 z + P_3 z^3,
\end{equation*}
where
\begin{align*}
  P_0 &= I, \\
  P_1 = P_{-1}^T &= \frac{1}{64}
  \begin{pmatrix}
    4 \sqrt{15} + 17 & & 4 \sqrt{15} + 16 \\
    -4 \sqrt{15} - 16 & &
    -4 \sqrt{15} - 17
  \end{pmatrix}, \\
  P_2 &= 0, \\
  P_3 = P_{-3}^T &= \frac{1}{64}
  \begin{pmatrix}
    15 - 4 \sqrt{15} & & 4 \sqrt{15} - 16 \\
    16 - 4 \sqrt{15} & & 4 \sqrt{15} - 15
  \end{pmatrix}.
\end{align*}
Numerical values of $P_k$ are given in table \ref{tab:02}.

\begin{table}
  \caption{The matrix coefficients $P_k$ of the half-band product filter
    $P$}
  \label{tab:02}
  \begin{tabular}{|c|rr|}
    \hline
    $k$ & \multicolumn{2}{c|}{$P_k$} \\
    \hline
    0 & 1 & \multicolumn{1}{l|}{0} \\
      & 0 & \multicolumn{1}{l|}{1} \\
    \hline
    1 
    &  $0.507686459137964$ &  $0.492061459137964$ \\
    & $-0.492061459137964$ & $-0.507686459137964$ \\
    \hline
    3 
    & $ -0.007686459137964$ & $ -0.007938540862036$ \\
    &  $ 0.007938540862036$ &  $0.007686459137964$ \\
    \hline
  \end{tabular}
\end{table}

\begin{figure}
  \includegraphics[width=4.5in]{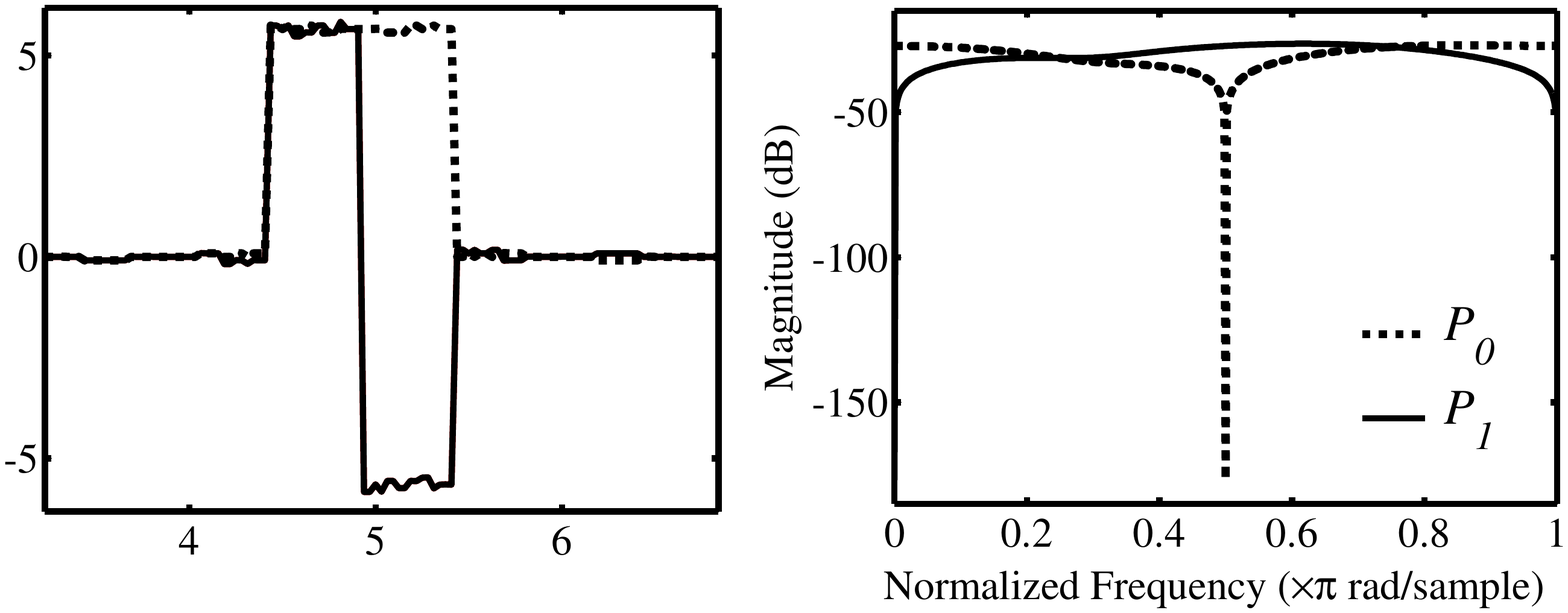} \\
  \caption{The impulse response (left) and frequency response
    (right) of the two components of the product filter $P = \left[
      P_0; P_1 \right]^T$.}
  \label{fig:04}
\end{figure}

The impulse responses of the product filter (see fig.~\ref{fig:04})
have the character of nearly Haar-type scaling and wavelet
functions. The frequency responses have good selectivity.

{\bf Step 3:} The matrix $T^{(f)}$ in this case looks like this:
\begin{equation*}
  T^{(f)} =
  \begin{pmatrix}
    P_0    & P_1^T & 0  & P_3^T & \\
    P_1 & P_0   & P_1^T   & 0   & P_3^T  & \\
    0 & & \ddots & & & \ddots &  \\
    P_3 &  & & \ddots & &  & P_3^T \\
    & \ddots & & & \ddots & & 0 \\
    &  & \ddots & &  & \ddots & P_1^T \\
    &  &  & P_3 & 0 & P_1 & P_0
  \end{pmatrix}
\end{equation*}

For any fixed $f$, the coefficients from the last row of
$L^{(f)}$ are approximate factors of $P(z)$. Suppressing the
dependence on $f$ in the notation, we define the filter
\begin{equation*}
  L^{SF}(z) = \sum_{k=0}^n C_k^{(f)} z^{-k}
\end{equation*}
and the matrix
\begin{align}\label{eq:LSF}
  L^{SF} &=
  \begin{bmatrix}
    C_0^{(f)} & C_1^{(f)} & C_2^{(f)} & C_3^{(f)} 
  \end{bmatrix} \\[5pt]
  &=
  \begin{bmatrix}
    L_{00} & L_{01} & L_{02} & L_{03} & L_{04} & L_{05} & L_{06} & L_{07} \\
    L_{10} & L_{11} & L_{12} & L_{13} & L_{14} & L_{15} & L_{16} & L_{17}
  \end{bmatrix}.
\end{align}
This notation will be used again later.

To measure the quality of the factorization, we compute
\begin{equation*}
  \Delta P(z) = P(z) - L^{SF}(z) \ L^{SF*}(z) = \sum_{k=-n}^n \Delta P_k \ z^k.
\end{equation*}
The residual is then
\begin{equation*}
  \Delta P = \max_{k=-n,\ldots,n} \max_{ij} | (\Delta P_k)_{ij} |.
\end{equation*}

In numerical experiments, the minimal error $\Delta P$ is
achieved for $f = 81$; the coefficients $C_k^{(81)}$ are given in
table \ref{tab:03}. The numerical error $\Delta P^{(81)} = 3.169 \cdot
10^{-9}$ for Bauer's method is much better than that of Wilson's
method for MSF, which cannot be lower then $10^{-5}$
\cite{ESS-17}. The impulse response of the obtained spectral factor is
non-regular, but the frequency response has good selectivity. This
approximate factor does not quite provide perfect reconstruction. This
is addressed in more detail in subsection \ref{subsec:inaccuracy}.

\begin{table}
  \caption{Coefficients of the lowpass spectral factor $L^{(f)}(z)$
    for $f = 81$, with $\Delta P^{(81)} = 3.169 \cdot 10^{-9}$.}
  \label{tab:03}
  \begin{tabular}{|c|rr|rr|}
    \hline
    $k$ & \multicolumn{2}{c|}{$C_k^{(81)}$} \\
    \hline
    $0$ 
    &  0.094428373297668 &  0 \\
    & -0.091754813647953 &  0.022310801334572 \\
    \hline
    $1$
    &  0.175943193428843 &  0.679731045265413 \\
    &  0.010360133828403 & -0.702056243347588 \\
    \hline
    $2$
    &  0.000129414745433 &  0.700756678587017 \\
    &  0.165030422671601 &  0.681046913905671 \\
    \hline
    $3$
    & -0.081190837390599 &  0.021002135513698 \\
    & -0.083853536287580 & -0.001275235049147 \\
    \hline
  \end{tabular}
\end{table}

Steps 4 and 5 will be covered in chapter \ref{sec:construction}.

\section{Construction of Approximate and Exact SA4 Multiwavelets}\label{sec:construction}

Above, we showed how to obtain the approximate spectral factor
$L^{SF}$ for a given size $f$. Given any factor of $P$, other possible
factors can be found be postmultiplication with suitable orthogonal
matrices, or in some cases by other manipulations.

In subsection \ref{subsec:approximate} we use a simple averaging
process and column reversal to produce a lowpass filter that is
similar to, but not identical to, the original SA4 multiscaling
function. We call this the {\em approximate SA4 multiwavelet}.

In subsection \ref{subsec:exact}, we add a rotation postfactor and an
averaging process to produce another multiwavelet which is even closer
to the original SA4 multiwavelet. We call this the {\em exact SA4
  multiwavelet}.

The newly derived multiscaling filters are denoted by
\begin{equation*}
  H^{(1)} = \left[ C_0^{(1)}, C_1^{(1)}, C_2^{(1)}, C_3^{(1)} \right]
\end{equation*}
for the
approximate SA4 multiwavelet, and
\begin{equation*}
    H^{(3)} = \left[ C_0^{(3)}, C_1^{(3)}, C_2^{(3)}, C_3^{(3)} \right]
\end{equation*}
for the exact SA4 multiwavelet. ($H^{(2)}$ is an intermediate step in
the calculation of $H^{(3)}$). In each case, we find the corresponding
multiwavelet filter by using a QR decomposition \cite{GvL-13}. The
multiscaling filter is factorized as
\begin{equation*}
  \begin{pmatrix}
    C_0 \\
    C_1 \\
    C_2 \\
    C_3
  \end{pmatrix} = QR =
  \begin{pmatrix}
    \vert   &        & \vert \\
    \boldq_1 & \cdots & \boldq_8 \\
    \vert   &        & \vert
  \end{pmatrix} \cdot
  \begin{pmatrix}
    I \\ 0 \\ 0 \\    0
  \end{pmatrix}
\end{equation*}
The coefficients $D_k$ are then obtained from the third and fourth
columns of $Q$.

We define the absolute error by
\begin{equation*}
  \Delta H^{(k)} = H - H^{(k)} = \left[ \Delta C_0^{(k)}, \Delta
    C_1^{(k)}, \Delta C_2^{(k)}, \Delta C_3^{(k)} \right],
\end{equation*}
where $\Delta C_n^{(k)} = C_n - C_n^{(k)}$ for $k = 1, 2, 3$, and
likewise for $\Delta G^{(k)}$.

To measure the deviation of the new filters from the original, we
introduce the {\em mean square error} (MSE) of the multiscaling and
multiwavelet functions,
\begin{equation*}
  \mbox{MSE-MF}^{(k)} = \sum_{ij} | \Delta H_{ij}^{(k)} |^2, \qquad
  \mbox{MSE-MwF}^{(k)} = \sum_{ij} | \Delta G_{ij}^{(k)} |^2,
\end{equation*}
and the{\em maximal absolute errors} (MAE) of the matrix coefficients,
the multiscaling function and the multiwavelet function as
\begin{align*}
  \mbox{MAE-MC}_\ell^{(k)} &= \max_{ij} | [\Delta C_\ell^{(k)}]_{ij} |, \\
  \mbox{MAE-MF} &= \max_{ij} | \Delta H_{ij}^{(k)} | = \max_{ij} | H_{ij} - H_{ij}^{(k)} |, \\
  \mbox{MAE-MwF} &= \max_{ij} | \Delta G_{ij}^{(k)} | = \max_{ij} | G_{ij} - G_{ij}^{(k)} |.
\end{align*}

\subsection{Approximate SA4 Multiwavelet}\label{subsec:approximate}

The spectral factor $L^{SF}$ obtained from Bauer's method leads to
non-symmetrical and non-regular functions. A simple algorithm can be
used to symmetrize and regularize the scaling functions. Below, we are
using the notation from eq.~\eqref{eq:LSF}. 

First, we average the absolute values of the first and fourth, as well
as the second and third, matrix coefficients of the spectral factor:
\begin{equation*}
   \begin{split}
     L_{0,6}^{SF} &= \frac{1}{4} \left( |L_{00}^{SF}| + |L_{10}^{SF}| + |L_{06}^{SF}| + |L_{16}^{SF}| \right), \\
     L_{1,7}^{SF} &= \frac{1}{4} \left( |L_{01}^{SF}| + |L_{11}^{SF}| + |L_{07}^{SF}| + |L_{17}^{SF}| \right), \\
     L_{2,4}^{SF} &= \frac{1}{4} \left( |L_{02}^{SF}| + |L_{12}^{SF}| + |L_{04}^{SF}| + |L_{14}^{SF}| \right), \\
     L_{3,5}^{SF} &= \frac{1}{4} \left( |L_{05}^{SF}| + |L_{15}^{SF}| + |L_{05}^{SF}| + |L_{15}^{SF}| \right).
   \end{split}
\end{equation*}
and construct average matrices
\begin{equation*}
     \begin{bmatrix}
       L_{2,4}^{SF} & L_{3,5}^{SF} \\
       L_{2,4}^{SF} & L_{3,5}^{SF}
     \end{bmatrix}, \qquad
     \begin{bmatrix}
       L_{0,6}^{SF} & L_{1,7}^{SF} \\
       L_{0,6}^{SF} & L_{1,7}^{SF}
     \end{bmatrix}.
\end{equation*}
Second, the matrix coefficients $C_1^{(1)}$ and $C_3^{(1)}$ are
obtained by multiplying the averaged matrices from the right by $J =
\diag(1,-1)$, but keep the signs, while the matrix coefficients
$C_0^{(1)}$ and $C_2^{(1)}$ are obtained by multiplying $C_1^{(1)}$
and $C_3^{(1)}$ on the left and right by $U = \mbox{antidiag}(1,1)$:
\begin{equation*}
   \begin{split}
     C_3^{(1)} &= \sign(L_3^{SF})
     \begin{bmatrix}
       L_{0,6}^{SF} & L_{1,7}^{SF} \\
       L_{0,6}^{SF} & L_{1,7}^{SF}
     \end{bmatrix} \cdot J, \\
     C_0^{(1)} &= U \cdot C_3^{(1)} \cdot U, \\
     C_1^{(1)} &= \sign(L_1^{SF})
     \begin{bmatrix}
       L_{2,4}^{SF} & L_{3,5}^{SF} \\
       L_{2,4}^{SF} & L_{3,5}^{SF}
     \end{bmatrix} \cdot J, \\
     C_2^{(1)} &= U \cdot C_1^{(1)} \cdot U.
   \end{split}
\end{equation*}
This produces coefficients quite close to the original coefficients
in~\eqref{eq:24}.

\begin{figure}
  \includegraphics[width=3in]{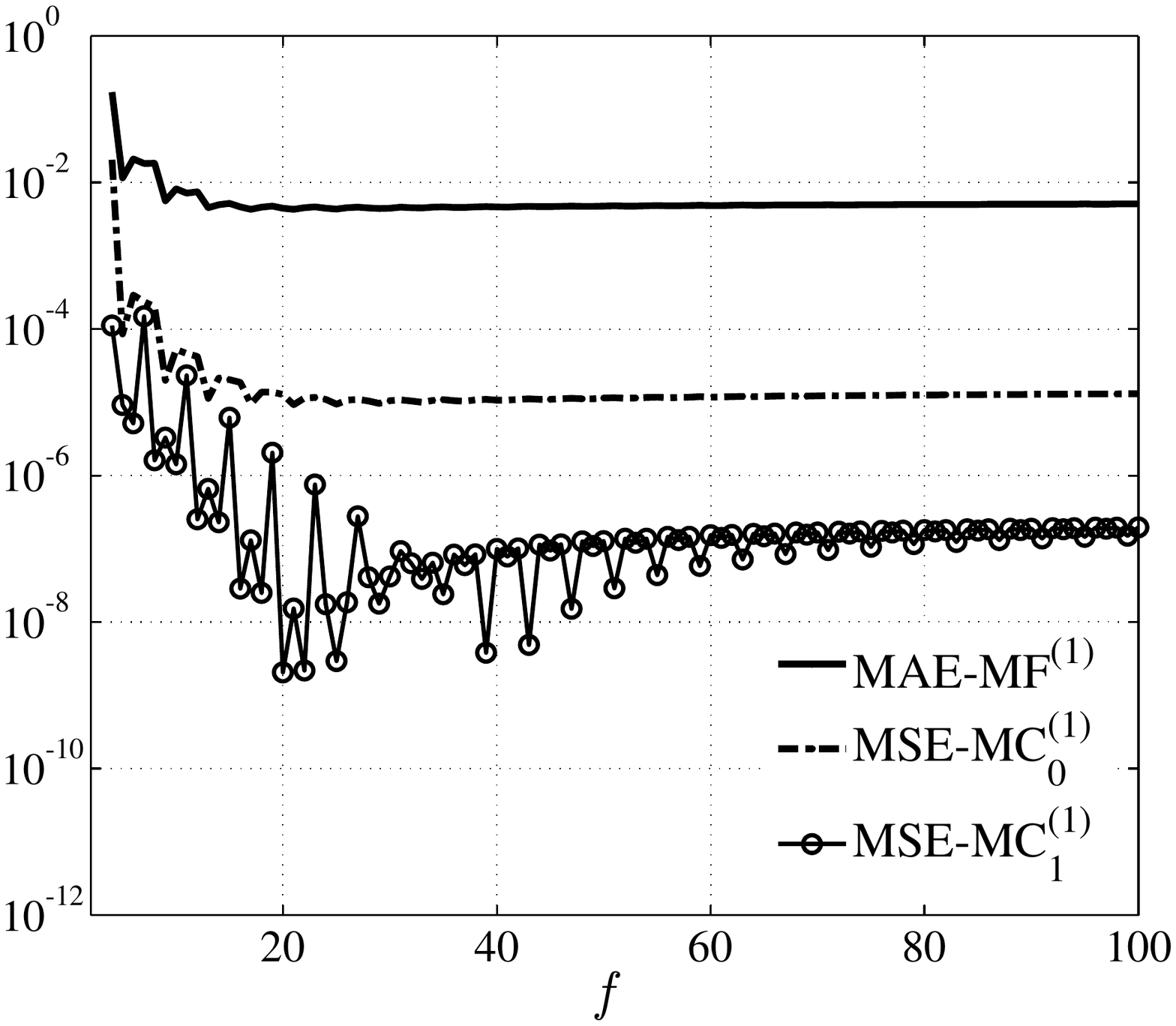}
  \caption{Dependence of maximum absolute error MAE-MF$^{(1)}$ and
    mean squared errors MSE-MC$_0^{(1)}$ = MSE-MC$_3^{(1)}$,
    MSE-MC$_1^{(1)}$ = MSE-MC$_2^{(1)}$ on matrix size $f$
    (log.~scale).}
  \label{fig:05}
\end{figure}

In order to determine the best approximate solution, we investigate
the influence of the size $f$ of the Toeplitz matrix $T^{(f)}$ for the
numerical errors MAE-MF$^{(1)}$, MSE-MC$_0^{(1)}$, and
MSE-MC$_1^{(1)}$. These errors are shown in fig.~\ref{fig:05} in
logarithmic scale, for $f$ up to 100. Note that MSE-MC$_3^{(1)}$ =
MSE-MC$_0^{(1)}$, and MSE-MC$_2^{(1)}$ = MSE-MC$_1^{(1)}$, by
construction.

The important minimal errors MAE-MF$^{(1)}$ and MSE-MwF$^{(1)}$ for
$f = 21$ are tabulated in table~\ref{tab:04}; the matrix coefficients
of are listed in table~\ref{tab:05}. Obviously, the development of
the errors shows the convergence of the algorithm, as well as the
dependency of the errors on the size $f$. It also shows that the main
error comes from the matrix coefficients $C_0^{(1)}$ and $C_3^{(1)}$,
rather than $C_1^{(1)}$ and $C_2^{(1)}$.

The error MSE-MF$^{(1)}$ behaves essentially the same way as
MAE-MC$_0^{(1)}$.  Unfortunately, the convergence of MAE-MC$_0^{(1)}$
and MAE-MC$_3^{(1)}$ stalls at a relatively large number. Thus, the
approximate SA4 multiwavelet is not very accurate. Better matrix
coefficients are achieved in the next subchapter.

\begin{table}
  \caption{The minimum MAEs and MSEs for the approximate and exact SA4
    multiscaling and multiwavelet functions, the size $f$ and the
    corresponding angle $\theta$.}
  \label{tab:04}

  \begin{tabular}{|cccccc|}
    \hline
    \hline
    $f$ & & MAE-MF$^{(1)}$ & MAE-MwF$^{(1)}$ & MSE-MF$^{(1)}$ & MSE-MwF$^{(1)}$ \\
    \hline
    21 & & \multicolumn{2}{c}{$4.2797 \cdot 10^{-3}$} & 
    \multicolumn{2}{c|}{$4.6139 \cdot 10^{-6}$} \\
    \hline
    \hline
    \multicolumn{1}{c}{ } \\
    \hline 
    \hline 
    $f$ & $\theta$  & MAE-MF$^{(2)}$ & MAE-MwF$^{(2)}$ & MSE-MF$^{(2)}$ & MSE-MwF$^{(2)}$ \\
    \hline
    15,400 & $1.444628609880879$ & \multicolumn{2}{c}{$1.220 \cdot 10^{-4}$}
    & \multicolumn{2}{c|}{$7.376 \cdot 10^{-9}$} \\
    \hline
    \hline 
    \multicolumn{1}{c}{ } \\
    \hline 
    \hline 
    $f$ & $\theta$  & MAE-MF$^{(3)}$ & MAE-MwF$^{(3)}$ & MSE-MF$^{(3)}$ & MSE-MwF$^{(3)}$ \\
    \hline
    12,042 & $1.444630025395427$ & $8.285 \cdot 10^{-8}$ & $8.557
    \cdot 10^{-8}$ & $3.835 \cdot 10^{-15}$ & $3.717 \cdot 10^{-15}$  \\    
    \hline
    \hline
  \end{tabular}
\end{table}

\begin{table}
  \caption{Coefficients of the approximate SA4 multiscaling and multiwavelet functions for $f = 21$.}
  \label{tab:05}
  \begin{tabular}{|c|rr|rr|}
    \hline
    n & \multicolumn{2}{c|}{$C_n^{(3)}$} & \multicolumn{2}{c|}{$D_n^{(3)}$} \\
    \hline
    $0$ 
    & $0.011165766264837$ &  $0.088552225597447$ & $-0.088552225597447$ &  $0.011165766264837$ \\
    & $0.011165766264837$ & $-0.088552225597447$ & $-0.088552225597447$ & $-0.011165766264837$ \\
    \hline
    $1$
    &  $0.691600252880066$ & $0.088718768987217$ &  $0.088718768987217$ & $-0.691600252880066$ \\
    & $-0.691600252880066$ & $0.088718768987217$ & $-0.088718768987217$ & $-0.691600252880066$ \\
    \hline
    $2$
    & $0.691600252880066$ & $-0.088718768987217$ & $0.088718768987217$ &  $0.691600252880066$ \\
    & $0.691600252880066$ &  $0.088718768987217$ & $0.088718768987217$ & $-0.691600252880066$ \\
    \hline
    $3$
    &  $0.011165766264837$ & $-0.088552225597447$ & $-0.088552225597447$ & $-0.011165766264837$ \\
    & $-0.011165766264837$ & $-0.088552225597447$ &  $0.088552225597447$ & $-0.011165766264837$ \\
    \hline
  \end{tabular}
\end{table}

\subsection{Exact SA4 Multiwavelet}\label{subsec:exact}

The second algorithm leads to more regular scaling functions. We
multiply the coefficients of $L^{SF}$ from the right by the
unitary matrix $U(\theta)$
\begin{equation*}
  C_k^{(2)} = C_k^{SF} \cdot U(\theta) = C_k^{SF} \cdot
  \begin{pmatrix}
    \cos\theta & \sin\theta \\
    -\sin\theta & \cos\theta
  \end{pmatrix}, \qquad k=0,\ldots,3.
\end{equation*}
This produces the approximate factor $H^{(2)}$. 

The angle $\theta$ (in radians) is calculated as the average value
\begin{equation*}
  \theta = \frac{1}{2} \left( \cos^{-1}(\mbox{even}) +
  \sin^{-1}(\mbox{odd}) \right),
\end{equation*}
where
\begin{equation*}
  \mbox{even} = \frac{1}{4} \left(
  \left| L_{00} + L_{04} + L_{10} + L_{14} \right| +
  \left| L_{02} + L_{06} + L_{12} + L_{16} \right| \right)
\end{equation*}
and
\begin{equation*}
  \mbox{odd} = \frac{1}{4} \left(
  \left| L_{01} + L_{05} + L_{11} + L_{15} \right| +
  \left| L_{03} + L_{07} + L_{13} + L_{17} \right| \right).
\end{equation*}
Again, we are using the notation from eq.~\eqref{eq:LSF}.

\begin{figure}
  \includegraphics[width=2.5in]{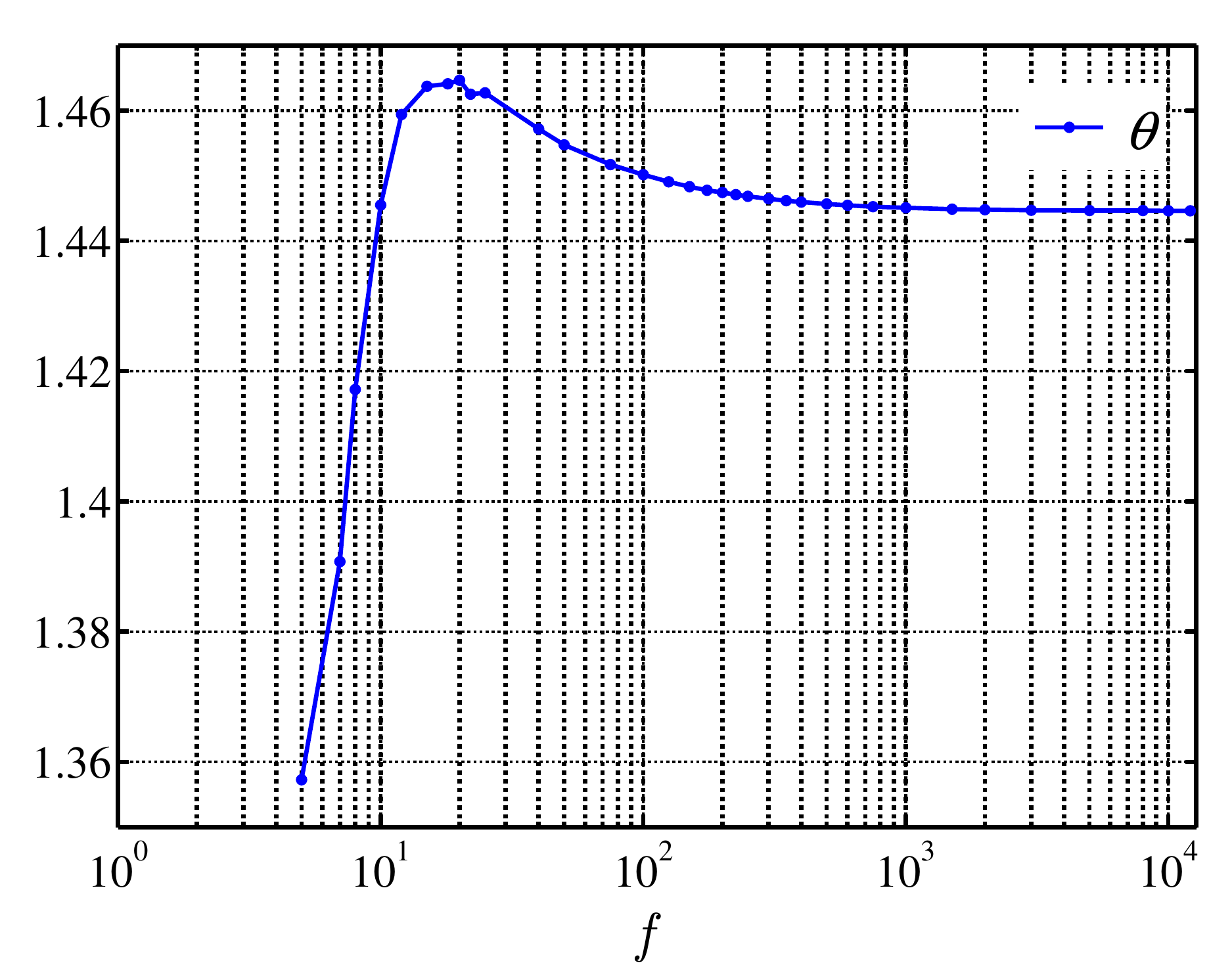}
  \caption{Dependence of the angle $\theta$ (radians) on size $f$ (log scale).}
  \label{fig:06}
\end{figure}

\begin{figure}
  \includegraphics[height=7in]{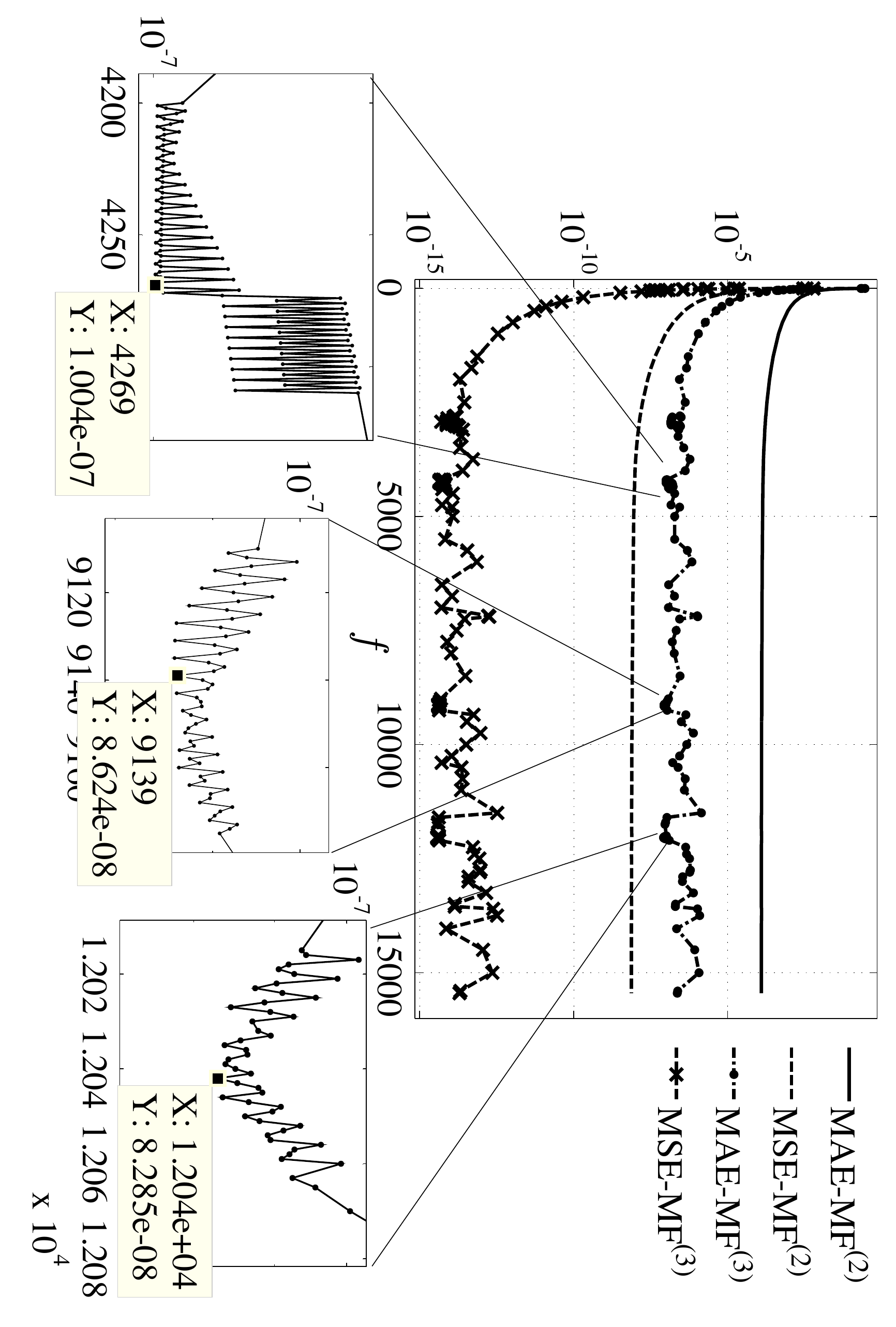}
  \caption{Dependence of the errors MSE-MF$^{(i)}$
    and MAE-MF$^{(i)}$ $i = 2, 3$, on the size of $f$ (log scale).}
  \label{fig:07}
\end{figure}

Figure~\ref{fig:06} shows the dependence of the angle $\theta$ on
the size $f$ (log.~scale). It can be seen that in the angle there is
one overshoot, after that it is decreasing very slowly. The minimum
errors are obtained at the maximal possible size of Toeplitz matrix $f
= 15,400$, with the values shown in table~\ref{tab:04}. This shows the
influence of the very slow convergence; despite right multiplication
with the unitary matrix, desirable precison cannot be an achieved. The
errors MSE-MF$^{(2)}$ and MAE-MF$^{(2)}$ for $H^{(2)}$ are shown in
fig.~\ref{fig:07}.

To increase the numerical precision and obtain the exact SA4
multiscaling function, we again apply an averaging approach.
\begin{align*}
  H^{(3)} &= [C_0^{(3)}, C_1^{(3)}, C_2^{(3)}, C_3^{(3)}], \\ 
  (C_0^{(3)})_{ij} &= \sign((C_0^{(2)})_{ij}) \cdot \frac{1}{2} \left( \left| (C_0^{(2)})_{ij} \right| + \left| (C_3^{(2)})_{ij} \right| \right), \\
  (C_1^{(3)})_{ij} &= \sign((C_1^{(2)})_{ij}) \cdot \frac{1}{2} \left( \left| (C_1^{(2)})_{ij} \right| + \left| (C_2^{(2)})_{ij} \right| \right), \\
  (C_2^{(3)})_{ij} &= \sign((C_2^{(2)})_{ij}) \cdot \frac{1}{2} \left( \left| (C_1^{(2)})_{ij} \right| + \left| (C_2^{(2)})_{ij} \right| \right), \\
  (C_3^{(3)})_{ij} &= \sign((C_3^{(2)})_{ij}) \cdot \frac{1}{2} \left( \left| (C_0^{(2)})_{ij} \right| + \left| (C_3^{(2)})_{ij} \right| \right).
\end{align*}

The errors for $H^{(3)}$ are shown in fig.~\ref{fig:07}, with minima
shown in table~\ref{tab:04}.  Let us consider the accuracy of the
multiscaling functions obtained at the leading local minima of the
error MAE-MF$^{(3)}$ . There are 3 minima, as can be seen in the
zoomed part of fig.~\ref{fig:07}. The first local minimum is at the
value $f = 4,269$, the second at $f = 9,139$, the third at
$f=12,042$. This means that increasing the size of the Toeplitz matrix
leads to slow improvement in the precision of the exact multiscaling
function. Nevertheless, due to applying the averaging method, the
minimal errors in $H^{(3)}$ and $G^{(3)}$ for $f = 12,042$ are smaller
by 2 or 3 orders of magnitude than the errors in $H^{(2)}$ and
$G^{(2)}$, and smaller by 5 orders of magnitude than the errors in
$H^{(1)}$ and $G^{(1)}$.

We choose $H^{(3)}$ for $f = 12,042$ for the coefficients of the exact SA4
multiwavelet, given in table \ref{tab:06}. They are almost the same as the
original coefficients in \eqref{eq:24}.

\begin{table}
  \scriptsize
  \caption{Coefficients of the exact SA4 multiscaling and multiwavelet functions for $f = 12,042$.}
  \label{tab:06}
  \begin{tabular}{|c|rr|rr|}
    \hline
    n & \multicolumn{2}{c|}{$C_n^{(3)}$} & \multicolumn{2}{c|}{$D_n^{(3)}$} \\
    \hline
    $0$ 
    & $0.01122679201336641$ &   $0.08838843031594616$ &  $-0.08838843285624247$ &   $0.01122679262734716$ \\
    & $0.01122679228175923$ &  $-0.08838843013543155$ &  $-0.08838843302944295$ &  $-0.01122679235802515$ \\
    \hline
    $1$
    &  $0.6958799459541121$ &   $0.08838843050191729$ &  $0.08838843089794699$ &  $-0.6958799676294173$ \\ 
    & $-0.6958799462085388$ &   $0.08838842817713602$ & $-0.08838843321541406$ &  $-0.6958799673174058$ \\
    \hline
    $2$
    & $0.6958799459541121$ &   $-0.08838843050191729$ &  $0.08838843089794699$ &    $0.6958799676294172$ \\
    & $0.6958799462085388$ &    $0.08838842817713602$ &  $0.08838843321541406$ &   $-0.6958799673174058$ \\
    \hline
    $3$
    &  $0.01122679201336641$ &  $-0.08838843031594616$ &  $-0.08838843285624255$ &  $-0.01122679262734722$ \\
    & $-0.01122679228175923$ &  $-0.08838843013543155$ &  $0.08838843302944295$ &  $-0.01122679235802536$ \\
    \hline
  \end{tabular}
\end{table}

\section{Performance Analysis}\label{sec:perf}

\subsection{Comparison With Other Multiwavelets}

Table~\ref{tab:07} gives the results of our experiments with the
approximate and exact SA4 multiwavelets, compared with other
well-known orthogonal and biorthogonal filter banks. Comparisons
include coding gain (CG), Sobolev smoothness (S) \cite{CMP-98},
symmetry/antisymmetry, and length. All of the multiwavelets considered
have symmetry/antisymmetry.

The CG for orthogonal transforms is a good indication of the
performance in signal processing. It is the ratio of arithmetic and
geometric means of channel variances $\sigma_i^2$:
\begin{equation*}
  CG = \frac{\frac{1}{r} \sum_{i=1}^r \sigma_i^2}{\left( \prod_{i=1}^r
  \sigma_i^2 \right)^{1/r}}
\end{equation*}

Coding gain is one of the most important factors to be considered in
many applications. It is always greater than or equal to 1; greater
values are better. CG is equal to 1 if all the variances are equal,
which means that it is not possible to clearly distinguish between the
smooth and the detailed components of the multiwavelet transformation
coefficients.  To estimate CG, the variance is computed using a first
order Markov model AR(1) with intersample autocorrelation coefficient
$\rho = 0.95$ \cite{DS-98}.

The Sobolev exponent $S$ of a filter bank measures the
$L^2$-differentiability of the corresponding multiscaling function
$\boldphi = [\phi_0, \phi_1]^T$, and thus also the multiwavelet
function $\boldpsi = [\psi_0, \psi_1]^T$. It is completely determined
by the multiscaling symbol $H(z)$.

The obtained Sobolev regularity and CGs of the approximate and exact
SA4 multiwavelets are equal: $S^{(1)} = S^{(3)} = 0.9919$, and
$CG^{(1)} = CG^{(3)} = 3.7323$ dB. This means that for some
applications we can use the SA4 multiwavelet.

According to table~\ref{tab:07}, the approximate and exact SA4
multiwavelets are better than most commonly used filter banks. They
also allow an economical lifting scheme for future implementation.

\begin{figure}
  \includegraphics[height=1.3in]{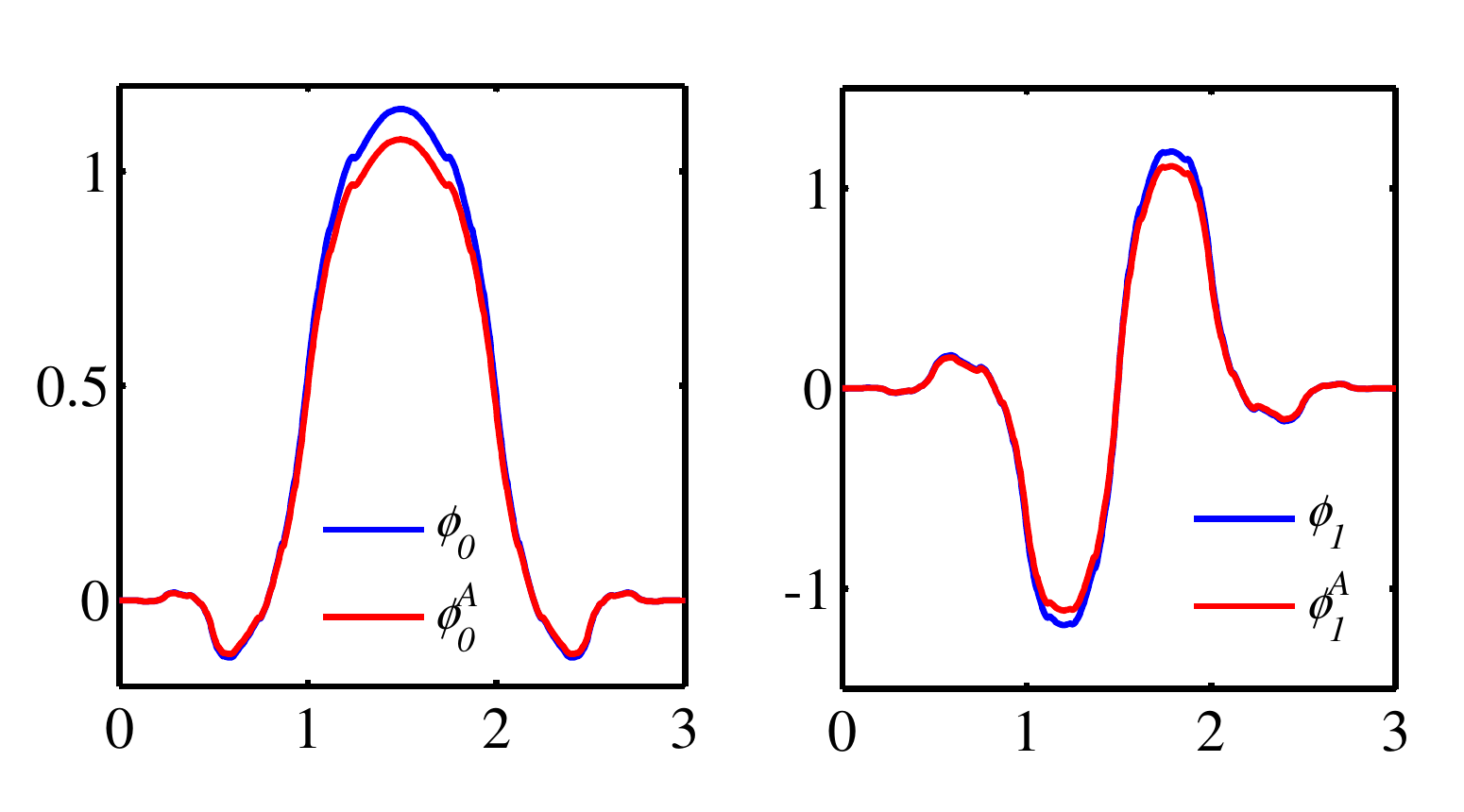}
  \includegraphics[height=1.3in]{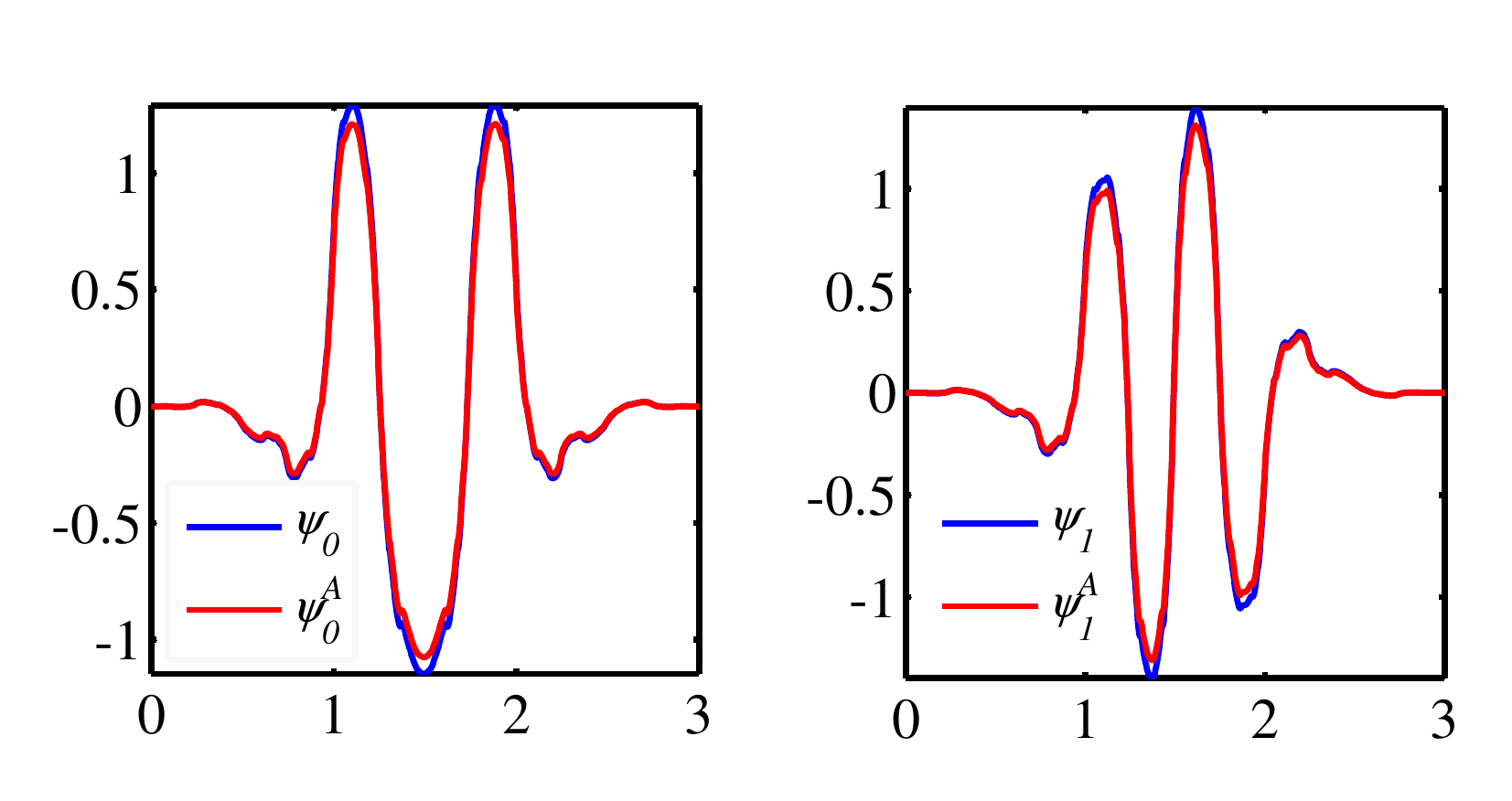} \\
  \begin{center}
  \includegraphics[height=1.6in]{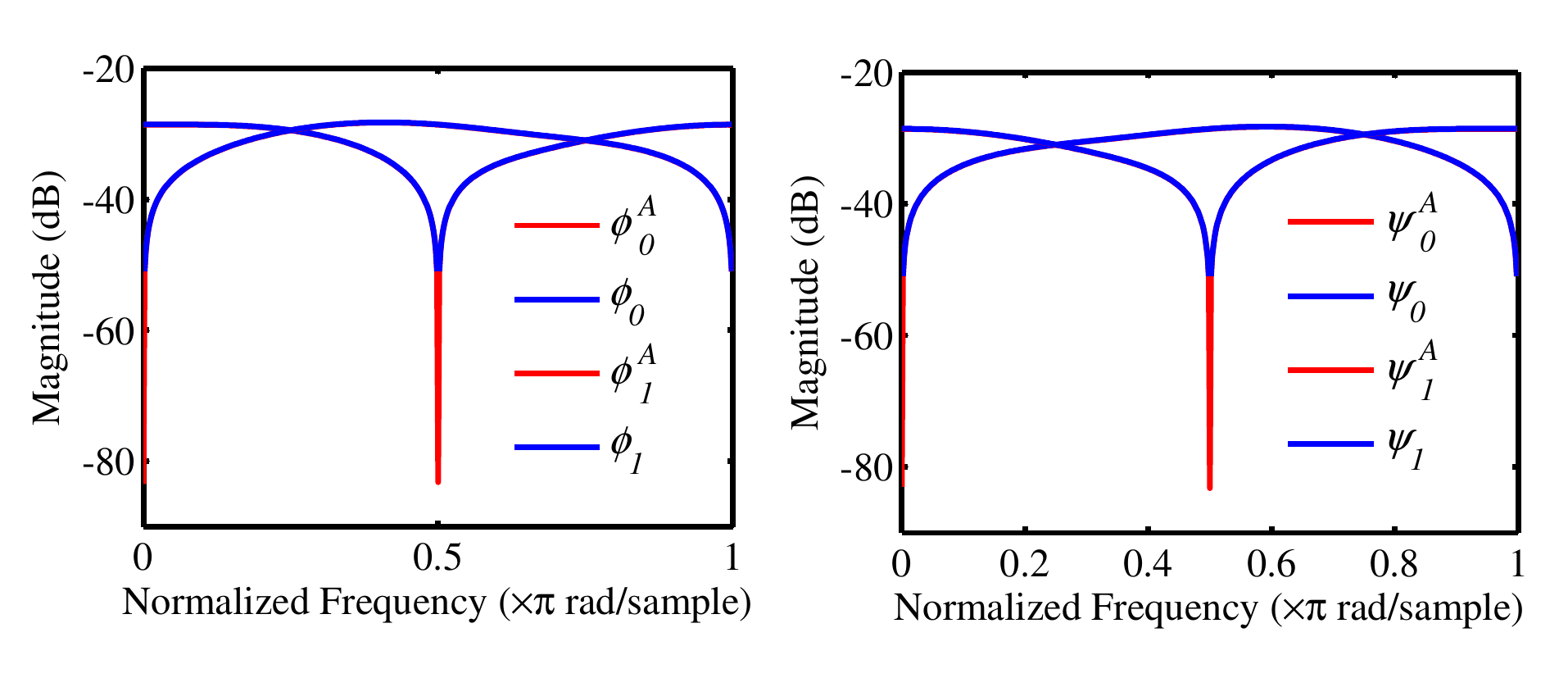} 
  \end{center}
  \caption{Top: The impulse responses of the components of the exact
    SA4 multiwavelet $\boldphi = (\phi_0, \phi_1)^T$, $\boldpsi =
    (\psi_0, \psi_1)^T$ {\em (in blue)}, compared with the approximate
    multiwavelet $\boldphi^A = (\phi_0^A, \phi_1^A)^T$, $\boldpsi^A =
    (\psi_0^A, \psi_1^A)^T$ {\em (in red)}; left to right: $\phi_0$,
    $\phi_1$, $\psi_0$, $\psi_1$; Bottom: The frequency response of
    the components of the exact SA4 multiwavelet {\em (in blue)}
    compared with the approximate multiwavelet {\em (in
      red)}; left: $\phi_0$, $\phi_1$; right: $\psi_0$, $\psi_1$.}
  \label{fig:08}
\end{figure}

Fig.~\ref{fig:08} shows a comparison of the approximate and exact SA4
multiwavelets in time domain and impulse response.  The influence of
the lower precision of the approximate SA4 multiwavelet can be seen in
the time domain image.

\begin{table}
  \caption{Comparison of coding gain (CG), Sobolev regularity (S),
    symmetry/antisymmetry (S/A), and length, for various
    multiwavelets.}
  \label{tab:07}
  \begin{tabular}{|c|c|c|c|c|}
    \hline
    Multifilter & CG & S & S/A & Length \\
    \hline
    \hline
    Biorthogonal Hermitian cubic spline (bih34n, \cite{S-98}) & 1.51 &	2.5 & yes & 3/5 \\
    \hline
    Biorthogonal (dual to Hermitian cubic spline) (bih32s, \cite{AZC-07}) & 1.01 & 2.5 & yes & 3/5 \\
    \hline
    Integer Haar \cite{CP-01} &	1.83 & 0.5 & yes & 2 \\
    \hline
    Chui-Lian (CL, \cite{CL-96}) & 2.06 & 1.06 & yes & 3 \\
    \hline
    Biorthogonal (dual to Hermitian cubic spline) (bih54n, \cite{AZC-07}) & 2.42 & 0.61 & yes & 5/3 \\
    \hline
    Biorthogonal (from GHM) (bighm2, \cite{AZC-07}) & 2.43 & 0.5 & yes & 2/6 \\
    \hline
    Biorthogonal (from GHM, dual to bighm2) (bighm6, \cite{AZC-07}) & 3.53 & 2.5 & yes & 6/2 \\
    \hline
    Biorthogonal (dual to Hermitian cubic spline) (bih52s, \cite{T-98}) & 3.69 & 0.83 & yes & 5/3 \\
    \hline
    Approximate SA4 & 3.73 & 0.99 & yes & 4 \\
    \hline
    Exact SA4 (\cite{STT-00}) & 3.73 & 0.99 & yes & 4 \\
    \hline
    Geronimo-Hardin-Massopust (GHM, \cite{GHM-94}) & 4.41 & 1.5 & yes & 4 \\
    \hline
    \end{tabular}
\end{table}

The magnitude of the approximate SA4 in the time domain is observably
smaller than for the exact SA4, while the frequency responses are
essentially identical (see fig.~\ref{fig:08}). Therefore, in applications where the
frequency response is important, we can use the approximate
multiscaling function.

\subsection{Influence of Inaccuracy of Filter Coefficients with Approximate SA4 Multiwavelet}
\label{subsec:inaccuracy}

The approximate SA4 multiwavelet $H^{(1)}$ is not quite an orthogonal
perfect reconstruction filter. In this section, we explore whether it
can still be useful in applications.

\subsubsection{1D Signal (No Additional Processing)}

The influence of the inaccuracy of the matrix filter coefficients in
the case of the approximate SA4 multiwavelet is shown by 1D
applications with no additional processing. By ``no additional
processing'' we mean that we simply decompose and reconstruct a signal
through several levels. No compression, denoising, or other processing
is done.

We use Haar balancing pre- and postfilters~\cite{CMP-98,SHSTH-99}
\begin{equation*}
  Q = \frac{1}{\sqrt{2}}
  \begin{pmatrix}
    1 & 1 \\
    -1 & 1
  \end{pmatrix},
\end{equation*}
whose small length preserves the time localization of the
multiwavelet decomposition, simplicity, orthogonality, and symmetry.

\begin{figure}
  (a) \includegraphics[width=4in]{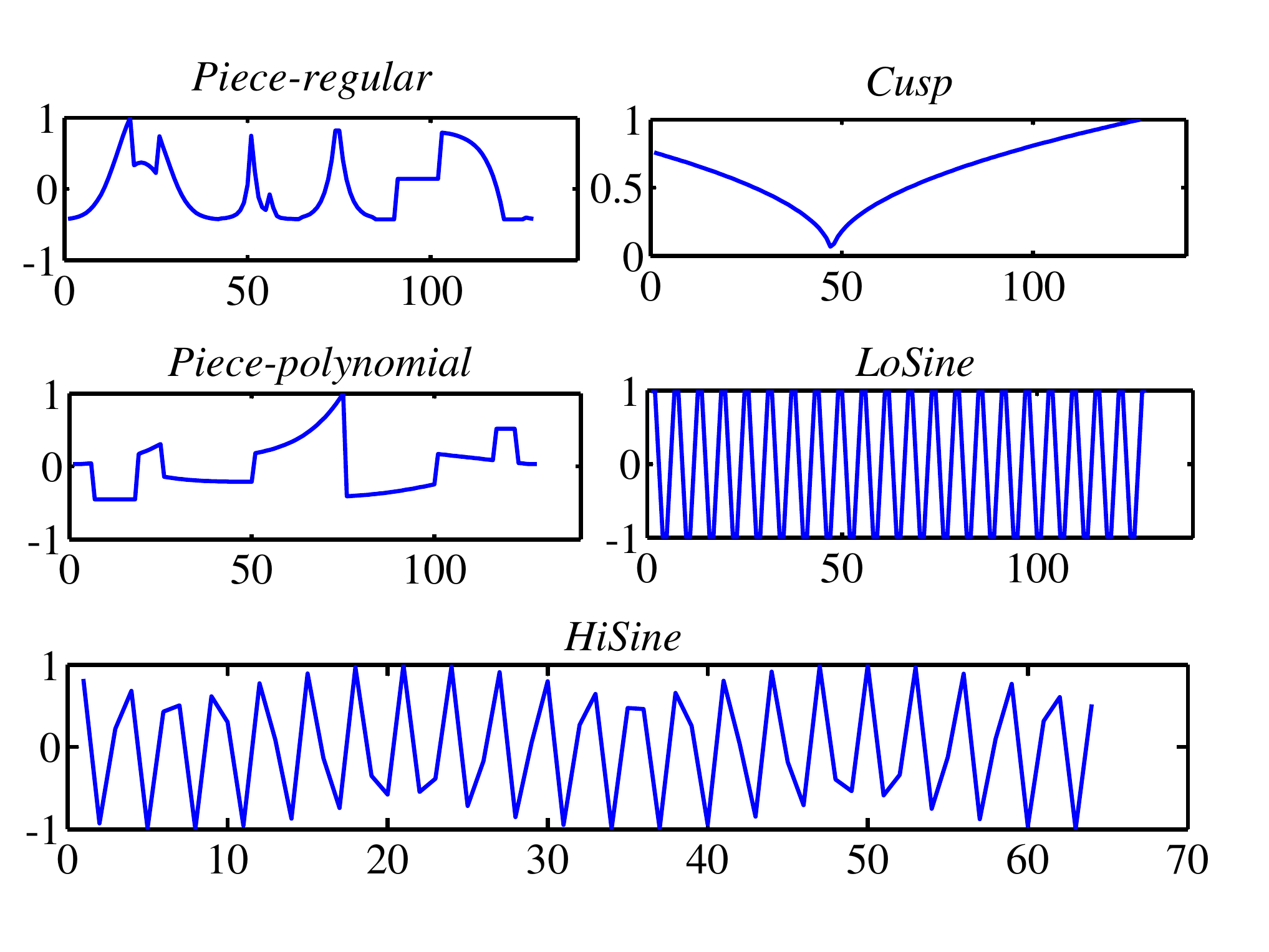} \\
  (b) \includegraphics[width=4.5in]{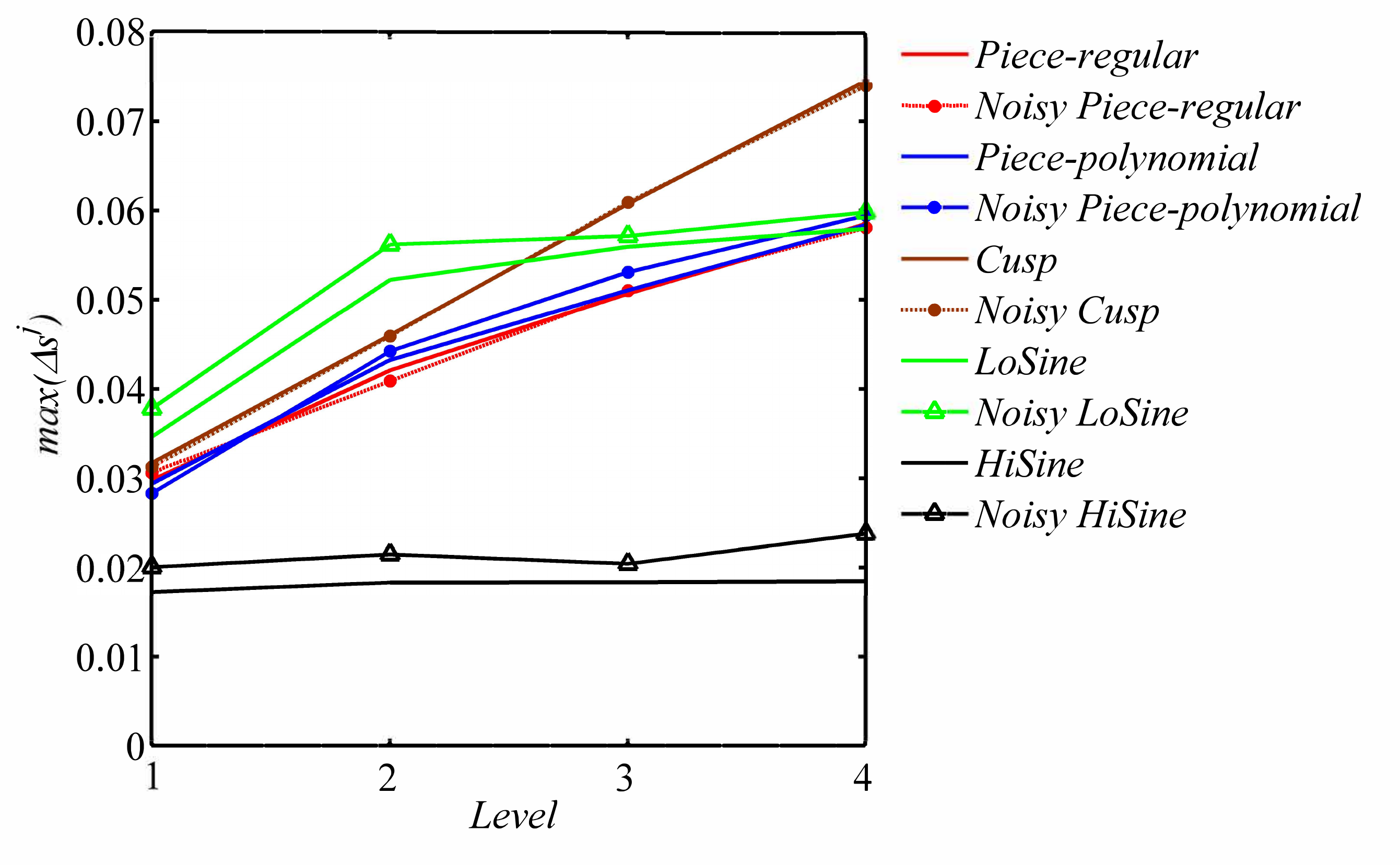}
  \caption{Decomposition and reconstruction with the balanced
    approximate SA4 multiwavelet; (a) The normed test signals; (b) The
    MAEs obtained at level $j$.}
  \label{fig:09}
\end{figure}

Results for the balanced version are shown in fig.~\ref{fig:09}(b).

The error in orthogonality of the approximate balanced SA4
multiwavelet is
\begin{equation*}
  \Delta H = I - Q \cdot H^{(1)} \cdot (H^{(1)})^T Q^T = 0.0117 \ I.
\end{equation*}

For the quality measure, we decompose and reconstruct five normed test
signals of length $N=2^7$, $s$ (without noise) and $\hat s = s +
\epsilon$ (with noise). The number of decomposition levels is 6.

The noise components $\epsilon_i$ are independent identically
distributed random variables with mean 0 and standard deviation
$\sigma$. We use maximum absolute errors (MAE)
\begin{equation*}
  \Delta s = \max_i \left| \mbox{reconstructed } s_i - \mbox{original  } s_i \right|,
\end{equation*}
where $s_i$ refers to the $i^{\mbox{\scriptsize th}}$ test signal.

The test signals are {\em 'Cusp'}, {\em 'HiSine'}, {\em 'LoSine'},
{\em 'Piece-regular'}, and {\em 'Piece-polynomial'}, implemented in
the Matlab environment. See fig.~\ref{fig:09}(a).

For the test signal {\em 'Cusp'}, increasing the number of
decomposition and reconstruction levels $j$ leads to a linear increase
of MAE. From the obtained small differences (about $3.5 \cdot
10^{-4}$) in error between noisy and noiseless signal, it follows that
the presence of noise does not make much difference in signals of this
type.

For the test signal {\em 'HiSine'} we observe nearly constant MAEs for
the noisy and noiseless signals with increasing level $j$; the
difference between the noisy and noiseless signal is in the interval
$[2-5.3] \cdot 10^{-3}$. Again, the presence of noise has only a
minimal influence.

For the test signal {\em 'LoSine'}, both noisy and noiseless MAEs show
a linear increase up to the second level, with only a small increase
at the third and fourth level of $-5.7 \cdot 10^{-3}$ and $3.6 \cdot
10^{-3}$, as shown in fig.~\ref{fig:09}(b). Therefore, after the
second level the influence of noise is weak.

For the test signal {\em 'Piece-regular'}, the MAE grows quadratically
with increasing $j$. The difference between noisy and noiseless
signals is smallest at a pre-/post-processing step ($6 \cdot 10^{-5}$
and largest at the second level ($1.2 \cdot 10^{-3}$). For the test
signal {\em 'Piece-polynomial'}, the MAE grows non-uniformly and
non-linearly with increasing $j$. The difference between noisy and
noiseless signals is smallest at the pre-/post-processing step ($3.8
\cdot 10^{-4}$) and largest at the third level ($2 \cdot 10^{-3}$).

\subsubsection{2D Signal (No Additional Processing)}

\begin{figure}
  (a) \includegraphics[width=1.2in]{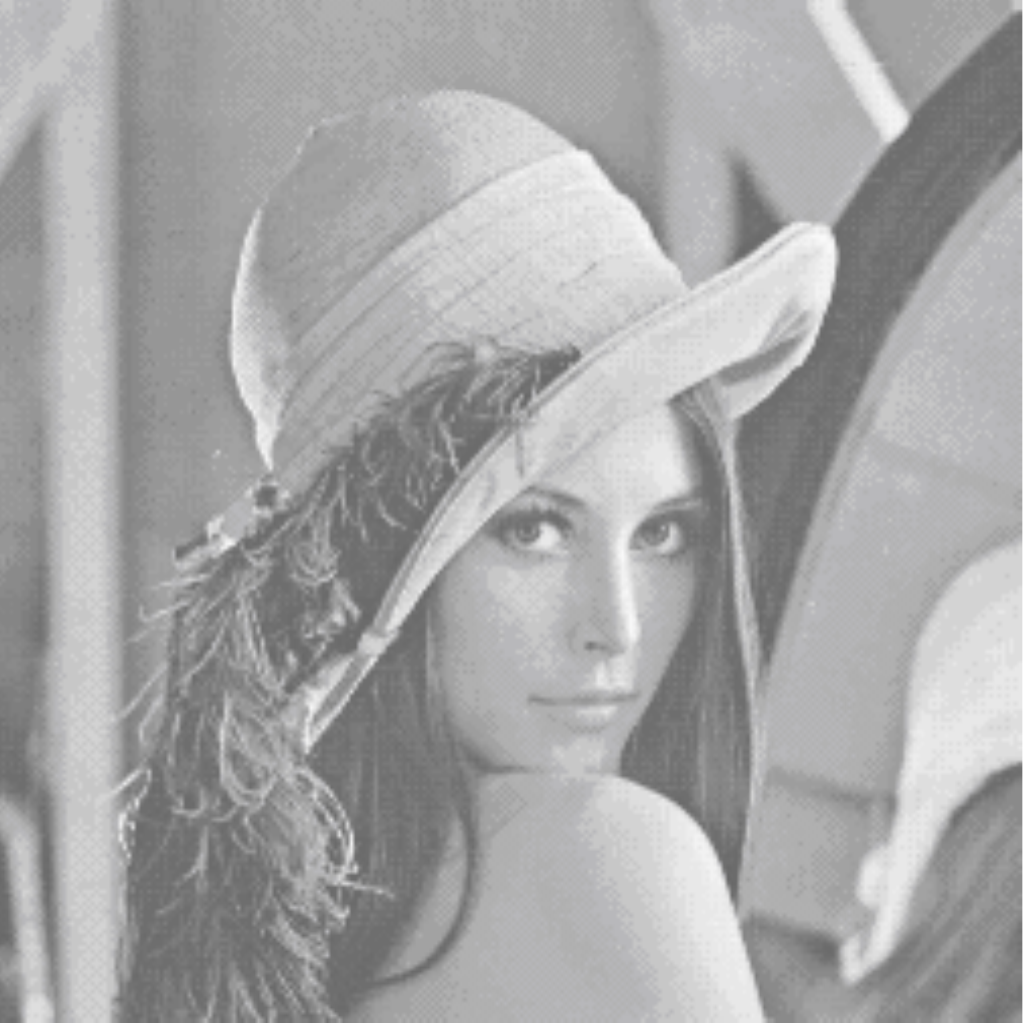}
  (b) \includegraphics[width=1.2in]{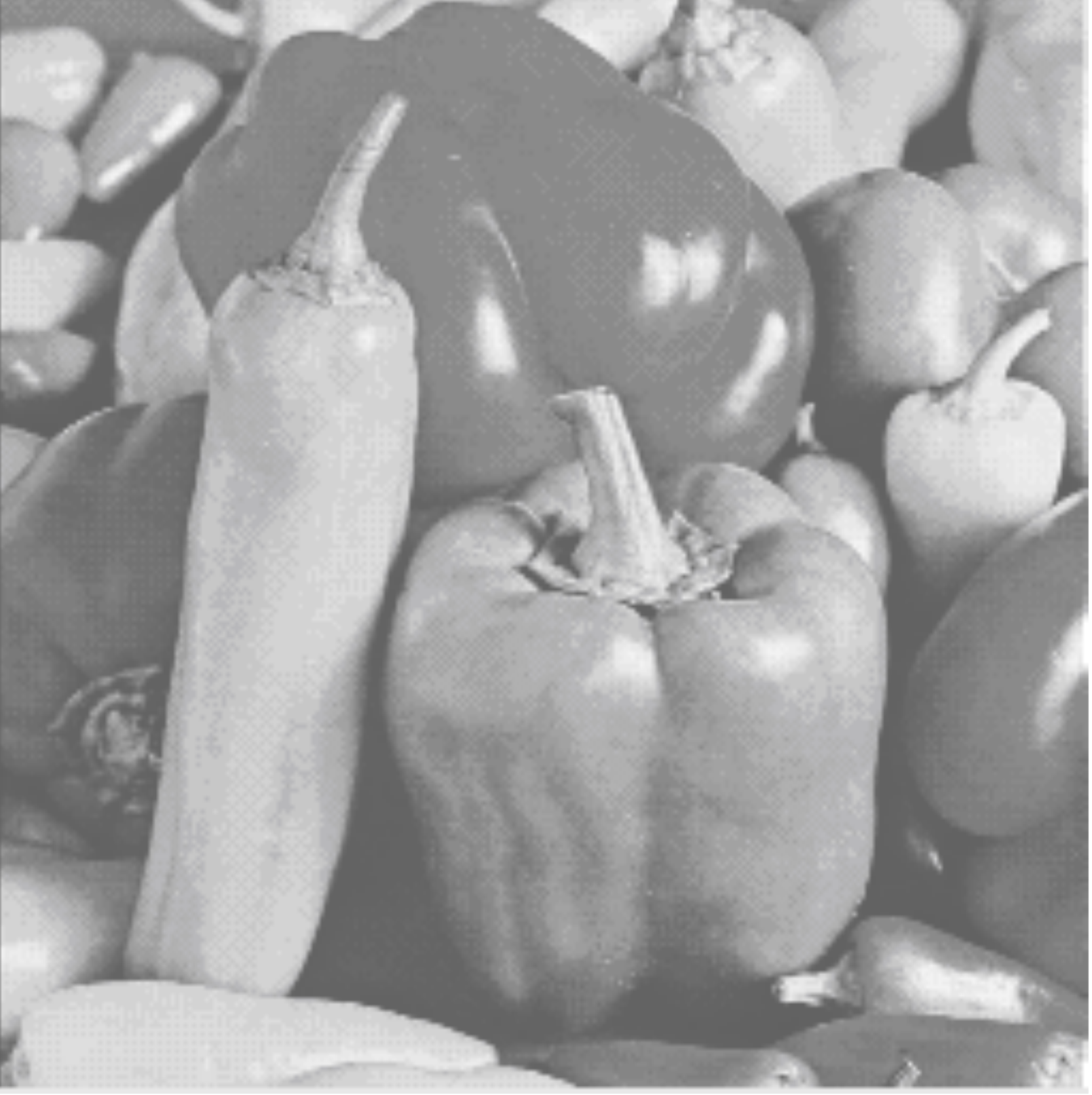}
  (c) \includegraphics[width=1.2in]{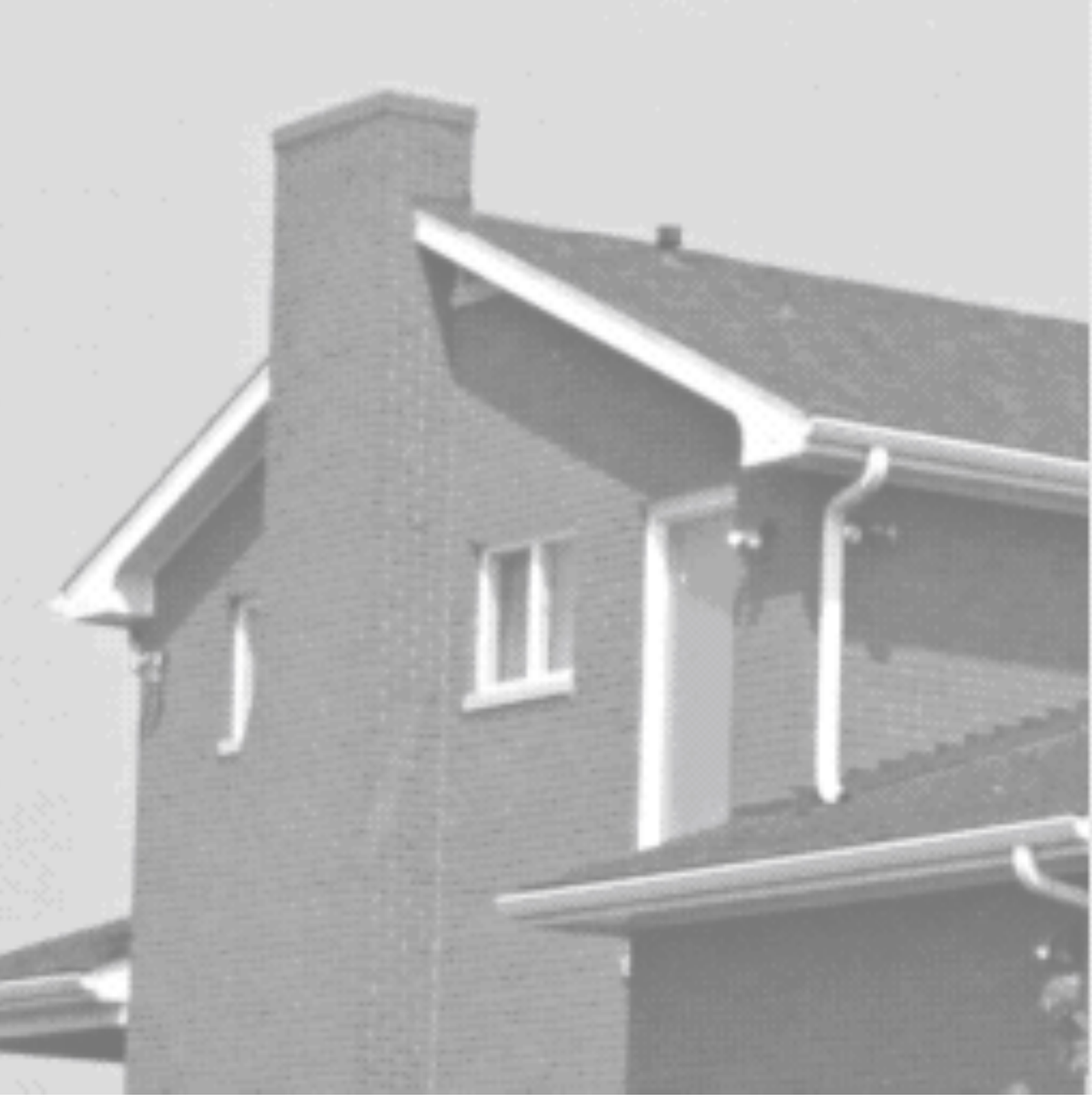} \\
  (d) \includegraphics[width=1.2in]{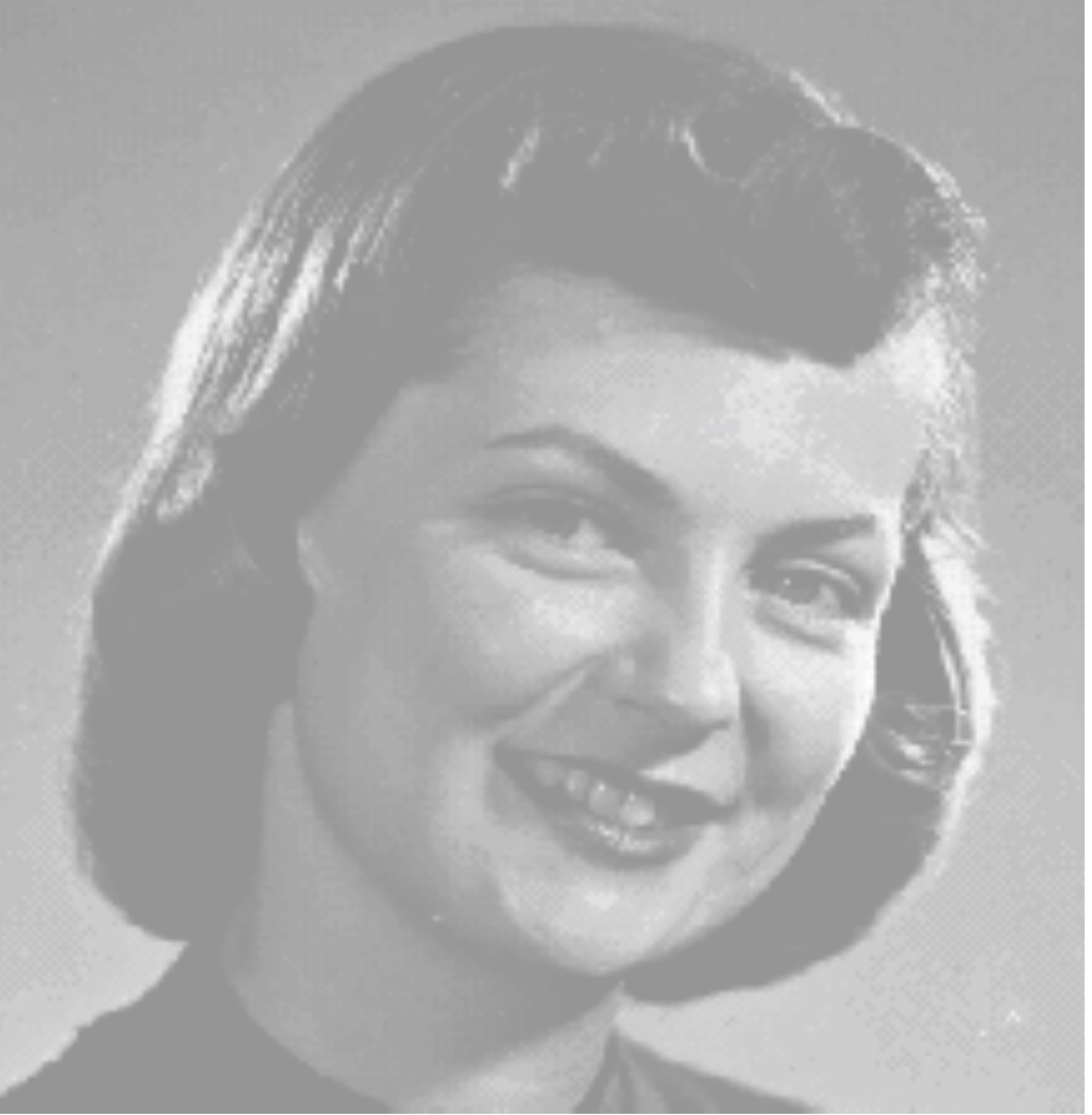}
  (e) \includegraphics[width=1.2in]{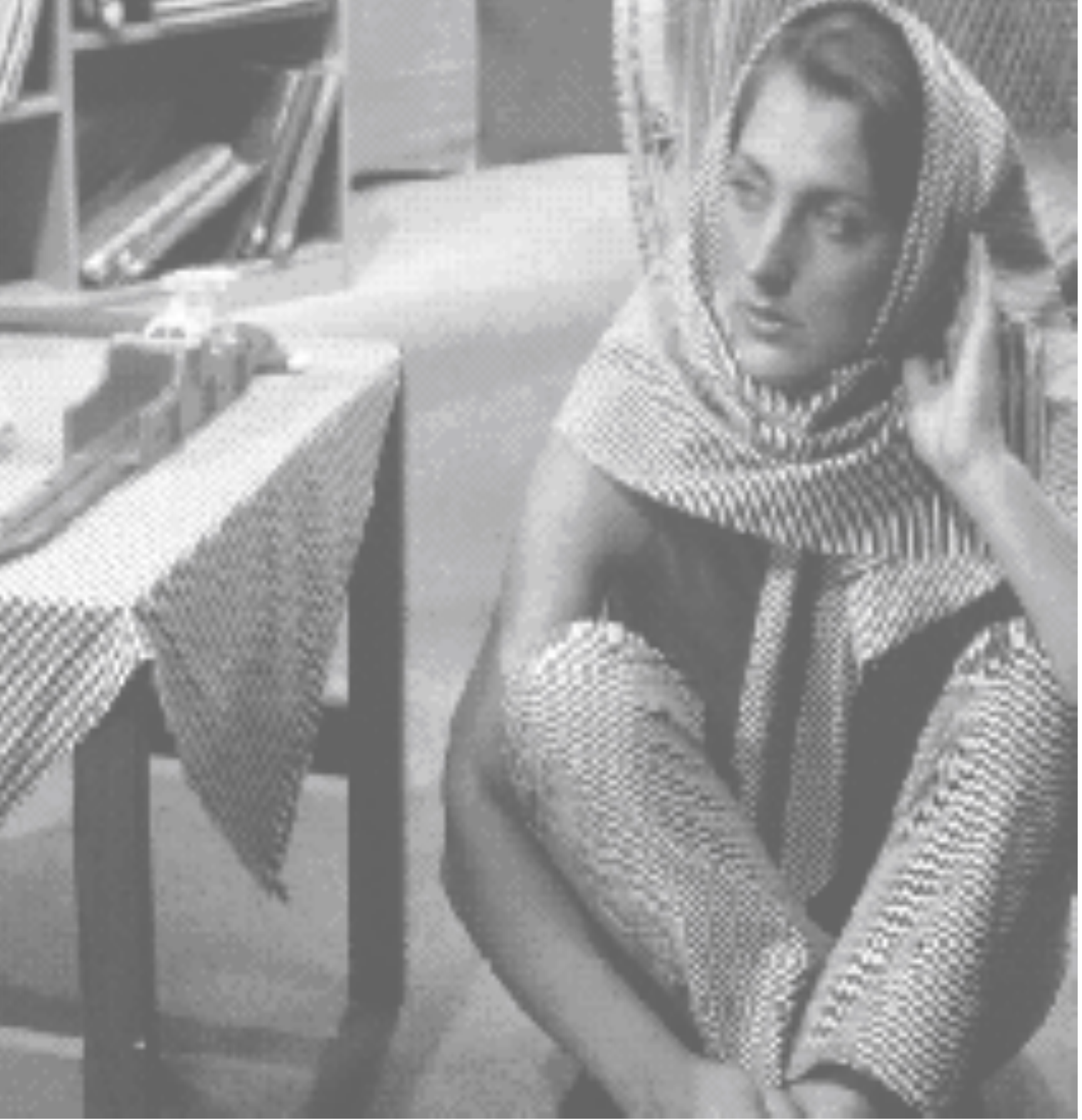}
  (f) \includegraphics[width=1.2in]{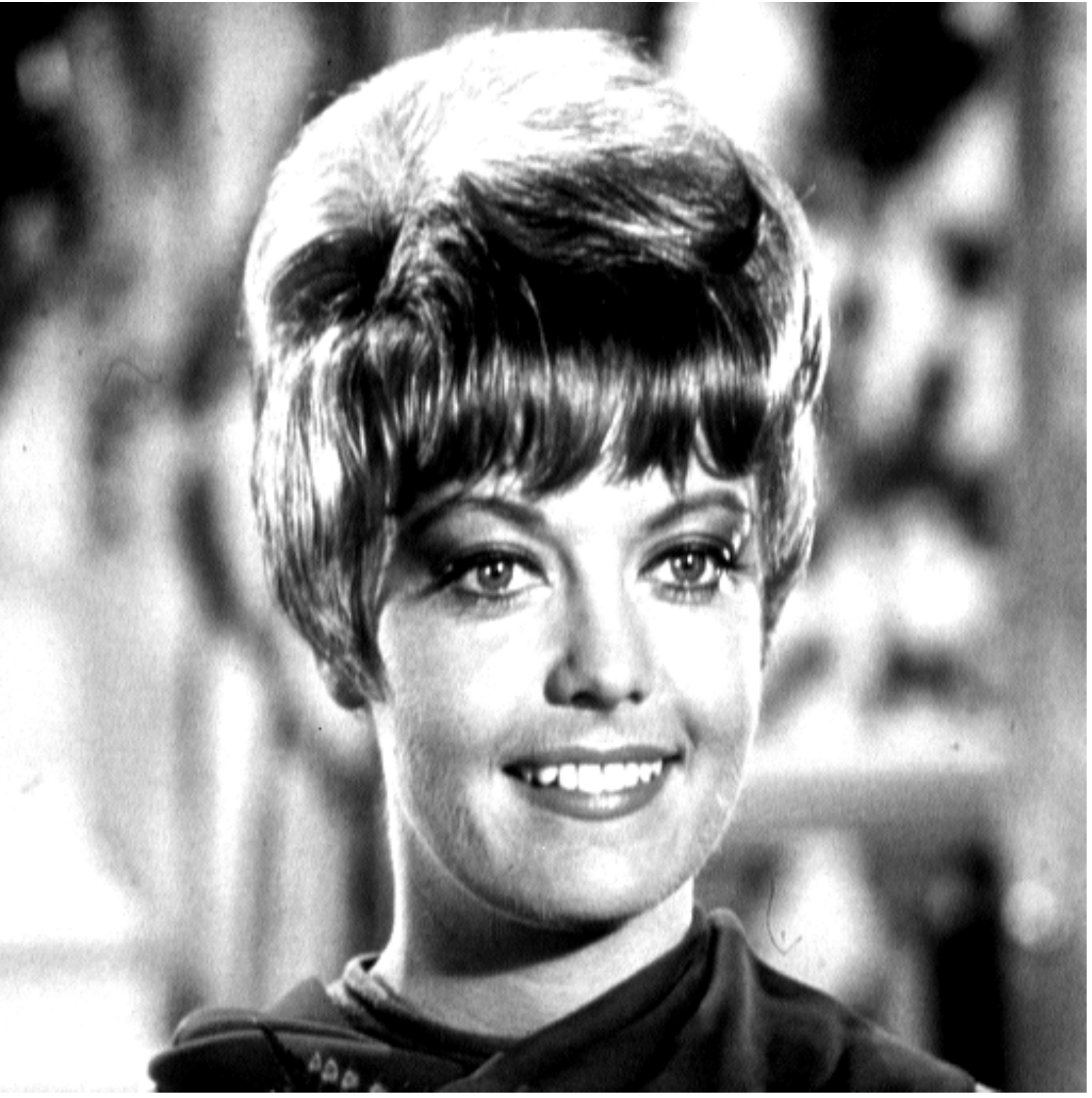}
  \caption{Test images of size $256 \times 256$ pixels) (a) Lena, (b)
    Peppers, and ($512 \times 512$ pixels) (c) House, (d) Girlface,
    (e) Barbara, and (f) Zelda.}
  \label{fig:10}
\end{figure}

The performance of the approximate SA4 multiwavelet was tested by
decomposing and reconstructing several gray level images, through 6 or
7 levels. Fig.~\ref{fig:10} shows the details for four of the images,
with the quality measure $PSNR = \log_{10}(255^2/MSE)$ dB.  The images
used are {\em 'Lena'} , {\em 'Peppers'}, {\em 'Girlface'}, and {\em
  'Barbara'}.

The approximate balanced SA4 multiwavelet applied to the four images
leads to an exponential decrease of the PSNRs. According to these
results, it is preferable to use no more than 3 or 4 levels of
decomposition.  After that, the reconstructed image has very low PSNRs
and visibly worse quality. In some applications, such as big data
archives, higher levels of decomposition may be useful.

\begin{figure}
  \includegraphics[width=3in]{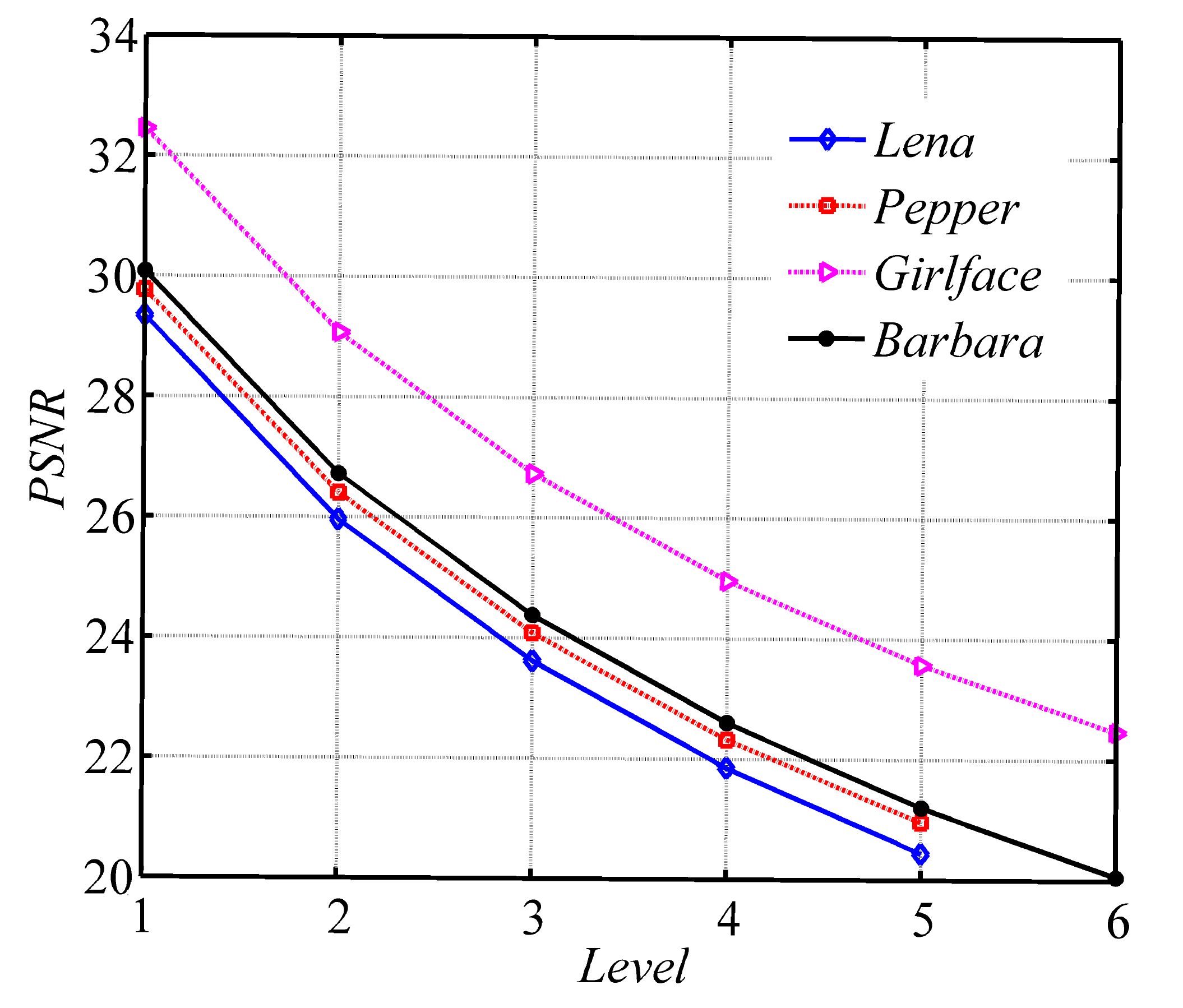}
  \caption{The PSNRs of four images for 2D decomposition and
    reconstruction without additional processing of 1 through 6 or 7
    levels with the approximate SA4 multiwavelet.}
  \label{fig:11}
\end{figure}

\subsection{Image Denoising}

In this section, we compare the balanced and non-balanced version of
the exact SA4 multiwavelet with the GHM multiwavelet (Downie and
Silverman 1998) and the Chui-Lian multiwavelet CL (Downie and
Silverman 1998), by considering image denoising with vector hard
thresholding (Donoho and Johnstone 1994), using 1-5 levels. The images
are {\em 'Lena'}, {\em 'Zelda'}, and {\em 'House'}, of size $512
\times 512$ pixels, with white additive Gaussian noise with variance
$\sigma=10$. See fig.~10.

The multiwavelet coefficients of the white noise are reduced at each
level, but uniformly distributed within each level.  Therefore, the
best approach to denoising is to find an appropriate threshold value
at each level.

\begin{figure}
  \includegraphics[height=2in]{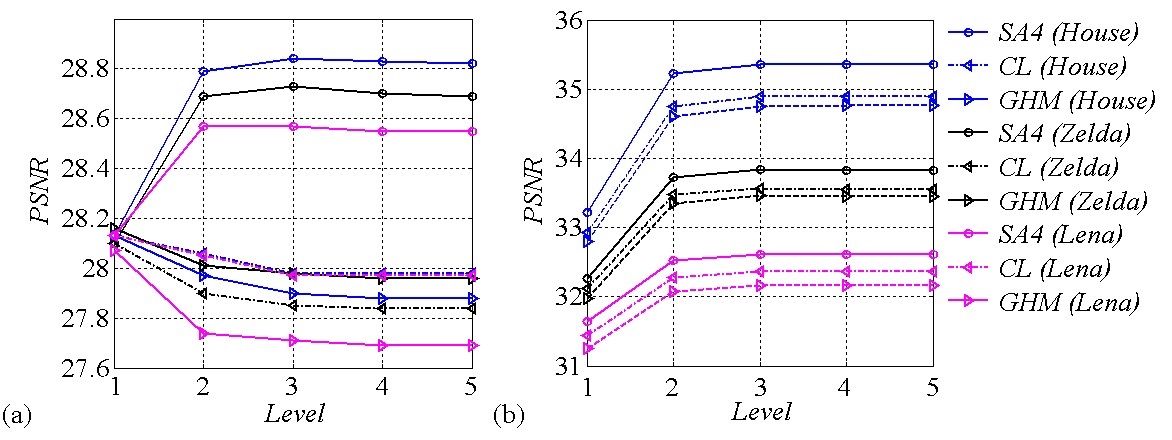} 
  \caption{Comparative analysis of PSNRs for image denoising with the
    exact SA4 multiwavelet and vector hard threshold through 2-6
    levels of the images {\em 'House'}, {\em 'Zelda'}, and {\em
      'Lena'} with $512 \times 512$ pixels and white additive
    Gaussian noise with variance $\sigma = 10$. Left: Non-balanced
    multiwavelets; Right: Balanced multiwavelets.}
  \label{fig:12}
\end{figure}

The exact SA4 multiwavelet, both balanced and non-balanced, achieved
the maximal PSNRs. The PSNR differences for the both versions of the
test images are $3.57$ dB at the pre-/post-processing step and $4.06$
dB at the fifth level for {\em ``Lena''}; $4.17$ dB
(pre-/post-processing step) and $5.13$ dB (fifth level) for {\em
  ``Zelda''}, $5.09$ dB (pre-/post-processing step) and $6.53$ dB
(fifth level) for {\em House''}.  Although balancing of multiwavelets
destroys the symmetry, it leads to increasing PSNRs and better image
denoising with the exact SA4 multiwavelet for the three test images
for the both version multiwavelets, while for the non-balanced GHM and
CL multiwavelets PSNRs decrease (see fig.~\ref{fig:12}(a)).

According to PSNRs, image decomposition and reconstruction through two
levels with the approximate SA4 multiwavelet is comparable to image
denoising with the non-balanced SA4 multiwavelet through three levels,
while one level is comparable with the balanced multiwavelet
pre-/post-processing step (see fig.~\ref{fig:11} and
fig.~\ref{fig:12}).

Obviously, the PSNRs of the exact balanced SA4 multiwavelet are
better than the well-known orthogonal multiwavelets GHM and
CL. They are useful in many applications.

\section{Conclusion}\label{sec:concl}

In this paper, we consider for the first time the problem of obtaining
an orthogonal multiscaling function by matrix spectral factorization
from a degenerate polynomial matrix. We show benchmark testing of MSF,
and apply Bauer's method to factoring the product filter of the SA4
multiwavelet. In addition, we show how to remove numerical errors and
improve the properties of the factors obtained from Cholesky
factorization, leading to fast convergent algorithms.

A very important part is obtaining the key angle $\theta$ in explicit
form. Based on the proposed averaging approach, we develop two filter
banks, the approximate and the exact SA4 orthogonal multiwavelets.

Experimental results have shown that the performances of the resulting
multiwavelets are better than those of the Chui-Lian multiwavelet and
biorthogonal multifilters, and are highly comparable to that of longer
multiwavelets. Theoretical analyses for the influence of the size $f$
of the Toeplitz matrix are considered, as well as simple 1D and 2D
applications.

After comparing both types of multifilters, we concluded that the
proposed averaging approach is a better way to remove numerical errors
and find the exact SA4 multiwavelet filter bank. It is important to
note that the performance of the balanced exact SA4 multiwavelet for
image denoising is better than the well-known orthogonal multiwavelets
GHM and CL, which are longer.

\begin{acknowledgements}
  The authors would like to thank the three anonymous referees for
  their critical review and helpful suggestions that allowed improving
  the exposition of the manuscript.
\end{acknowledgements}

\bibliographystyle{spmpsci}      
\bibliography{sa4}   

\begin{thebibliography}{10}
\providecommand{\url}[1]{{#1}}
\providecommand{\urlprefix}{URL }
\expandafter\ifx\csname urlstyle\endcsname\relax
  \providecommand{\doi}[1]{DOI~\discretionary{}{}{}#1}\else
  \providecommand{\doi}{DOI~\discretionary{}{}{}\begingroup
  \urlstyle{rm}\Url}\fi

\bibitem{AZC-07}
Averbuch, A.Z., Zheludev, V.A., Cohen, T.: Multiwavelet frames in signal space
  originated from {H}ermite splines.
\newblock IEEE Trans. Signal Process. \textbf{55}(3), 797--808 (2007)

\bibitem{B-56}
Bauer, F.L.: Beitr\"age zur {E}ntwicklung numerischer {V}erfahren f\"ur
  programmgesteuerte {R}echenanlagen. {II}. {D}irekte {F}aktorisierung eines
  {P}olynoms.
\newblock Bayer. Akad. Wiss. Math.-Nat. Kl. S.-B. \textbf{1956}, 163--203
  (1957) (1956)

\bibitem{C-07}
Charoenlarpnopparut, C.: One-dimensional and multidimensional spectral
  factorization using gr\"obner basis approach.
\newblock In: 2007 Asia-Pacific Conference on Communications, pp. 201--204
  (2007)

\bibitem{CP-01}
Cheung, K.W., Po, L.M.: Integer multiwavelet transform for lossless image
  coding.
\newblock In: Intelligent Multimedia, Video and Speech Processing, 2001.
  Proceedings of 2001 International Symposium on, pp. 117--120 (2001)

\bibitem{CL-96}
Chui, C.K., Lian, J.A.: A study of orthonormal multi-wavelets.
\newblock Appl. Numer. Math. \textbf{20}(3), 273--298 (1996).
\newblock Selected keynote papers presented at 14th IMACS World Congress
  (Atlanta, GA, 1994)

\bibitem{CDP-97}
Cohen, A., Daubechies, I., Plonka, G.: Regularity of refinable function
  vectors.
\newblock J. Fourier Anal. Appl. \textbf{3}(3), 295--324 (1997)

\bibitem{CNKS-96}
Cooklev, T., Nishihara, A., Kato, M., Sablatash, M.: Two-channel multifilter
  banks and multiwavelets.
\newblock In: Conference Proceedings, IEEE International Conference on
  Acoustics, Speech, and Signal Processing ICASSP-96, vol.~5, pp. 2769--2772
  (1996)

\bibitem{CMP-98}
Cotronei, M., Montefusco, L.B., Puccio, L.: Multiwavelet analysis and signal
  processing.
\newblock IEEE Transactions on Circuits and Systems II: Analog and Digital
  Signal Processing \textbf{45}(8), 970--987 (1998)

\bibitem{DJ-94}
Donoho, D.L., Johnstone, I.M.: Ideal spatial adaptation by wavelet shrinkage.
\newblock Biometrika \textbf{81}(3), 425--455 (1994)

\bibitem{DGHM-96}
Donovan, G.C., Geronimo, J.S., Hardin, D.P., Massopust, P.R.: Construction of
  orthogonal wavelets using fractal interpolation functions.
\newblock SIAM J. Math. Anal. \textbf{27}(4), 1158--1192 (1996)

\bibitem{DS-98}
Downie, T.R., Silverman, B.W.: The discrete multiple wavelet transform and
  thresholding methods.
\newblock IEEE Transactions on Signal Processing \textbf{46}(9), 2558--2561
  (1998)

\bibitem{DXD-13}
Du, B., Xu, X., Dai, X.: Minimum-phase {FIR} precoder design for multicasting
  over {MIMO} frequency-selective channels.
\newblock Journal of Electronics (China) \textbf{30}(4), 319--327 (2013)

\bibitem{EJL-09}
Ephremidze, L., Janashia, G., Lagvilava, E.: A simple proof of the
  matrix-valued {F}ej{\'e}r-{R}iesz theorem.
\newblock Journal of Fourier Analysis and Applications \textbf{15}(1), 124--127
  (2009)

\bibitem{ESS-17}
Ephremidze, L., Saied, F., Spitkovsky, I.: On the algorithmization of
  {J}anashia-{L}agvilava matrix spectral factorization method.
\newblock Submitted to IEEE Trans. Inform. Theory

\bibitem{F-05}
Fischer, R.F.: Sorted spectral factorization of matrix polynomials in {MIMO}
  communications.
\newblock IEEE Transactions on Communications \textbf{53}(6), 945--951 (2005)

\bibitem{GM-05}
Gan, L., Ma, K.K.: On minimal lattice factorizations of symmetric-antisymmetric
  multifilterbanks.
\newblock IEEE Trans. Signal Process. \textbf{53}(2, part 1), 606--621 (2005)

\bibitem{GHM-94}
Geronimo, J.S., Hardin, D.P., Massopust, P.R.: Fractal functions and wavelet
  expansions based on several scaling functions.
\newblock J. Approx. Theory \textbf{78}(3), 373--401 (1994)

\bibitem{GvL-13}
Golub, G.H., Van~Loan, C.F.: Matrix computations, fourth edn.
\newblock Johns Hopkins Studies in the Mathematical Sciences. Johns Hopkins
  University Press, Baltimore, MD (2013)

\bibitem{HCW-10}
Hansen, M., Christensen, L.P.B., Winther, O.: Computing the minimum-phase
  filter using the {QL}-factorization.
\newblock IEEE Trans. Signal Process. \textbf{58}(6), 3195--3205 (2010)

\bibitem{HHS-04}
Hardin, D.P., Hogan, T.A., Sun, Q.: The matrix-valued {R}iesz lemma and local
  orthonormal bases in shift-invariant spaces.
\newblock Advances in Computational Mathematics \textbf{20}(4), 367--384 (2004)

\bibitem{HLSH-07}
Hsung, T.C., Lun, D.P.K., Shum, Y.H., Ho, K.C.: Generalized discrete
  multiwavelet transform with embedded orthogonal symmetric prefilter bank.
\newblock IEEE Trans. Signal Process. \textbf{55}(12), 5619--5629 (2007)

\bibitem{HSLS-03}
Hsung, T.C., Sun, M.C., Lun, D.K., Siu, W.C.: Symmetric prefilters for
  multiwavelets.
\newblock IEE Proceedings-Vision, Image and Signal Processing \textbf{150}(1),
  59--68 (2003)

\bibitem{HM-12}
Huo, G., Miao, L.: Cycle-slip detection of {GPS} carrier phase with methodology
  of {SA4} multi-wavelet transform.
\newblock Chinese Journal of Aeronautics \textbf{25}(2), 227--235 (2012)

\bibitem{JK-85}
Je\v{z}ek, J., Ku\v{c}era, V.: Efficient algorithm for matrix spectral
  factorization.
\newblock Automatica \textbf{21}(6), 663--669 (1985)

\bibitem{J-98b}
Jiang, Q.: On the regularity of matrix refinable functions.
\newblock SIAM J. Math. Anal. \textbf{29}(5), 1157--1176 (1998)

\bibitem{J-98a}
Jiang, Q.: Orthogonal multiwavelets with optimum time-frequency resolution.
\newblock IEEE Trans. Signal Process. \textbf{46}(4), 830--844 (1998)

\bibitem{L-93}
Lawton, W.: Applications of complex valued wavelet transforms to subband
  decomposition.
\newblock IEEE Transactions on Signal Processing \textbf{41}(12), 3566--3568
  (1993)

\bibitem{LV-98}
Lebrun, J., Vetterli, M.: Balanced multiwavelets theory and design.
\newblock IEEE Trans. Signal Process. \textbf{46}(4), 1119--1125 (1998)

\bibitem{LV-01}
Lebrun, J., Vetterli, M.: High-order balanced multiwavelets: theory,
  factorization, and design.
\newblock IEEE Trans. Signal Process. \textbf{49}(9), 1918--1930 (2001)

\bibitem{LP-11}
Li, B., Peng, L.: Balanced multiwavelets with interpolatory property.
\newblock IEEE Trans. Image Process. \textbf{20}(5), 1450--1457 (2011)

\bibitem{LY-10}
Li, Y.F., Yang, S.Z.: Construction of paraunitary symmetric matrices and
  parametrization of symmetric orthogonal multiwavelet filter banks.
\newblock Acta Math. Sinica (Chin. Ser.) \textbf{53}(2), 279--290 (2010)

\bibitem{MRvF-96}
Massopust, P.R., Ruch, D.K., Van~Fleet, P.J.: On the support properties of
  scaling vectors.
\newblock Appl. Comput. Harmon. Anal. \textbf{3}(3), 229--238 (1996)

\bibitem{MS-97}
Micchelli, C.A., Sauer, T.: Regularity of multiwavelets.
\newblock Adv. Comput. Math. \textbf{7}(4), 455--545 (1997)

\bibitem{PS-98}
Plonka, G., Strela, V.: Construction of multiscaling functions with
  approximation and symmetry.
\newblock SIAM J. Math. Anal. \textbf{29}(2), 481--510 (1998)

\bibitem{LR-86}
Roux, J.L.: {2D} spectral factorization and stability test for {2D} matrix
  polynomials based on the radon projection.
\newblock In: Acoustics, Speech, and Signal Processing, IEEE International
  Conference on ICASSP '86., vol.~11, pp. 1041--1044 (1986)

\bibitem{SK-01}
Sayed, A.H., Kailath, T.: A survey of spectral factorization methods.
\newblock Numer. Linear Algebra Appl. \textbf{8}(6-7), 467--496 (2001).
\newblock Numerical linear algebra techniques for control and signal processing

\bibitem{STT-00}
Shen, L., Tan, H.H., Tham, J.Y.: Symmetric-antisymmetric orthonormal
  multiwavelets and related scalar wavelets.
\newblock Appl. Comput. Harmon. Anal. \textbf{8}(3), 258--279 (2000)

\bibitem{SB-86}
Smith, M., Barnwell, T.: Exact reconstruction techniques for tree-structured
  subband coders.
\newblock IEEE Transactions on Acoustics, Speech, and Signal Processing
  \textbf{34}(3), 434--441 (1986)

\bibitem{S-98}
Strela, V.: A note on construction of biorthogonal multi-scaling functions.
\newblock In: Wavelets, multiwavelets, and their applications ({S}an {D}iego,
  {CA}, 1997), \emph{Contemp. Math.}, vol. 216, pp. 149--157. Amer. Math. Soc.,
  Providence, RI (1998)

\bibitem{SHSTH-99}
Strela, V., Heller, P.N., Strang, G., Topiwala, P., Heil, C.: The application
  of multiwavelet filterbanks to image processing.
\newblock IEEE Transactions on Image Processing \textbf{8}(4), 548--563 (1999)

\bibitem{SW-00}
Strela, V., Walden, A.: Signal and image denoising via wavelet thresholding:
  orthogonal and biorthogonal, scalar and multiple wavelet transforms.
\newblock In: Nonlinear and nonstationary signal processing ({C}ambridge,
  1998), pp. 395--441. Cambridge Univ. Press, Cambridge (2000)

\bibitem{TST-99}
Tan, H.H., Shen, L.X., Tham, J.Y.: New biorthogonal multiwavelets for image
  compression.
\newblock Signal Processing \textbf{79}(1), 45--65 (1999)

\bibitem{TSLT-98}
Tham, J.Y., Shen, L., Lee, S.L., Tan, H.H.: Good multifilter properties: a new
  tool for understanding multiwavelets.
\newblock In: Proc. Intern. Conf. on Imaging, Science, Systems and Technology
  CISST-98, Las Vegas, USA, pp. 52--59 (1998)

\bibitem{TSLT-00}
Tham, J.Y., Shen, L., Lee, S.L., Tan, H.H.: A general approach for analysis and
  application of discrete multiwavelet transforms.
\newblock IEEE Trans. Signal Process. \textbf{48}(2), 457--464 (2000)

\bibitem{T-98}
Turcajov\'a, R.: Hermite spline multiwavelets for image modeling.
\newblock In: Proc. SPIE 3391, Wavelet Applications V, 46, Orlando, FL, USA,
  pp. 46--56 (1998)

\bibitem{WMW-15}
Wang, Z., McWhirter, J.G., Weiss, S.: Multichannel spectral factorization
  algorithm using polynomial matrix eigenvalue decomposition.
\newblock In: 2015 49th Asilomar Conference on Signals, Systems and Computers,
  pp. 1714--1718 (2015)

\bibitem{WM-57}
Wiener, N., Masani, P.: The prediction theory of multivariate stochastic
  processes. {I}. {T}he regularity condition.
\newblock Acta Math. \textbf{98}, 111--150 (1957)

\bibitem{W-72}
Wilson, G.T.: The factorization of matricial spectral densities.
\newblock SIAM J. Appl. Math. \textbf{23}, 420--426 (1972)

\bibitem{WLXL-10}
Wu, G., Li, D., Xiao, H., Liu, Z.: The {$M$}-band cardinal orthogonal scaling
  function.
\newblock Appl. Math. Comput. \textbf{215}(9), 3271--3279 (2010)

\bibitem{YK-78}
Youla, D.C., Kazanjian, N.N.: Bauer-type factorization of positive matrices and
  the theory of matrix polynomials orthogonal on the unit circle.
\newblock IEEE Trans. Circuits and Systems \textbf{CAS-25}(2), 57--69 (1978)

\end{thebibliography}

\end{document}